\documentclass[review,12pt]{elsarticle}
\usepackage{epstopdf}
\usepackage{amsfonts}
\usepackage{amsmath}
\usepackage{psfrag}
\usepackage{graphicx}
\usepackage{amssymb}
\usepackage{subfigure}
\usepackage{algorithmic}
\usepackage{algorithm}
\usepackage{color}
\usepackage{cases}
\usepackage{epsfig}
\usepackage{enumitem}
\usepackage{rotating}
\usepackage{amsmath,amssymb}
\usepackage{latexsym,verbatim}
\usepackage{pifont}
\usepackage{multirow}
\usepackage{mathrsfs}

\usepackage{hyperref}

\usepackage{booktabs}
\usepackage{bm}
\usepackage{pifont}
\usepackage{rotating}

\usepackage[T1]{fontenc}
\renewcommand{\vec}[1]{\boldsymbol{#1}}
\newcommand{\cmark}{\ding{51}}%
\newcommand{\xmark}{\ding{55}}%

\newproof{pot}{Proof}

\begin{document}

\begin{frontmatter}

\title{Adaptive Zeroing-Type Neural Dynamics for Solving Quadratic Minimization and Applied to Target Tracking
\tnoteref{t1}}

\author{Huiting He$^{1}$\corref{}}
\author{Chengze Jiang$^{2}$\corref{}}
\author{Yudong Zhang$^{3}$\corref{cor1}}
\author{Xiuchun Xiao$^{1}$\corref{cor1}}
\author{Zhiyuan Song$^{1}$\corref{}}
\address{$^{1}$ College of Electronic and Information Engineering, Guangdong Ocean University, Zhanjiang 524088, China}

\address{$^{2}$ School of Cyber Science and Engineering, Southeast University 210096, China}
\address{$^{3}$ School of Computing and Mathematical Sciences, University of Leicester, Leicester, LE1 7RH, UK}

\cortext[cor1]{Corresponding author. \\{\it E-mail addresses:} yudongzhang@ieee.org;xiaoxc@gdou.edu.cn.}

\begin{abstract}
{The time-varying quadratic minimization (TVQM) problem, as a hotspot currently, urgently demands a more reliable and faster--solving model. To this end, a novel adaptive coefficient constructs framework is presented and realized to improve the performance of the solution model, leading to the adaptive zeroing-type neural dynamics (AZTND) model. Then the AZTND model is applied to solve the TVQM problem. The adaptive coefficients can adjust the step size of the model online so that the solution model converges faster. At the same time, the integration term develops to enhance the robustness of the model in a perturbed environment. Experiments demonstrate that the proposed model shows faster convergence and more reliable robustness than existing approaches. Finally, the AZTND model is applied in a target tracking scheme, proving the practicality of our proposed model.}
\end{abstract}

\begin{keyword}
{Time-varying quadratic minimization (TVQM)\sep Zeroing-type neural dynamics (ZTND)\sep Adaptive coefficient}
\end{keyword}
\end{frontmatter}

\section{Introduction}\label{sec.intro}
As classical mathematical problems, some applications in the realms of energy system design \cite{AppEnergy}, image processing \cite{NCZNNCompare}, and robot kinematics \cite{AppRobort}, are abstractly modeled as quadratic minimization (QM) problems. In conclusion, for general QM problems, the traditional approach solves them by employing some numerical or iterative algorithms \cite{VPZNNCompare}. A message-passing scheme for solving QM problems is presented in \cite{AppMessagePass} by Ruozzi and Tatikonda. In addition, Zhang $et~al.$ present a QM-based dual-arm cyclic-motion-generation manipulator control scheme and analyze its properties from the perspective of cybernetics \cite{AppRobort}. It is worth noting that although considerable research has been devoted to solving conventional QM problems, studies aimed explicitly at the time-varying quadratic minimization (TVQM) problem are insufficient. Traditional solutions have serious lag errors when facing large-scale time-varying issues, resulting in inadequate solution accuracy and even the collapse of the solution system \cite{XXCOne}.
\par

To break through the dilemma that traditional algorithms cannot effectively deal with time-varying problems, Zhang $et~al.$ design the original zeroing neural network (OZNN) model \cite{ZNNProposed}. The OZNN model employs the derivative information of the time-varying problem to predict its evolution direction and continuously adjust the solution strategy of the solution system through a named evolution function \cite{XJ}. Therefore, with the help of the OZNN model, it can cope with various time-varying problems with extremely high accuracy \cite{LiuMei}. The OZNN model has been successfully applied to signal processing and automatic control fields due to its high accuracy and real-time solution advantages \cite{Qi}. However, the chief drawback of the OZNN model is its sensitivity to measurement noise. Its solution accuracy will be reduced distinctly in the presence of noise interference \cite{Wei}. Besides, the scale parameter of the OZNN model requires manually set and tuned. Hence tedious and repeated adjustments are necessary when facing actual engineering application problems \cite{XJ}. In recent years, it has been reported that much research aims to optimize the ZNN models. A versatile recurrent neural network termed the VRNN model is presented by Xiao $et~al.$  to solve the TVQM problem \cite{AppWideUse}, conquers the drawbacks of the OZNN model, takes time to converge when it approaches infinity, and accelerates the model to globally converge within a finite time \cite{FTZNNXiao}. On this basis, theoretical analysis, including the predefined-time convergence of the strictly predefined-time convergent ZNN (PTCZNN) model, is performed with the convergence of models \cite{PTCZNNCompare}. Noteworthily, both the VRNN model and the PTCZNN model implement the convergence properties of the control solution model by constructing a special activation function and not synthesizing the scale parameter in models.  Different from these fixed-valued neural-dynamic models, a varying-parameter convergent-differential neural network  (VP-CDNN) model is based on time-variant incremental scale parameters \cite{VPZNNCompare}. The VP-CDNN model converges exponentially and maintains better robustness under perturbation conditions than the OZNN model. However, in model implementation and engineering applications, as the time continues to increase, the monotonically increasing time-varying design parameters may be too large to be achieved or violate the objective limit, which results in the solving failure \cite{Huang}. When implementing zeroing-type models, it is unavoidable that models may be interfered with various measurement noises, leading to the reduction of system solution accuracy and even the collapse of the solution system. To this end, Jin $et~al.$ present the modified ZNN (MZNN) model in \cite{MZNNCompare}, which introduces the integral information into the solution evolution formula for the first time. However, parameters of the MZNN model still require tedious manual adjustments, which leads to a lot of additional computational resources and redundant adjustment processes \cite{Alawad}. He $et~al.$ present a residual learning framework to simplify the training process of deep neural networks, which explicitly reformulates layers as learning residual functions concerning the layer inputs \cite{DeepRL}. Inspired by the residual learning framework and combined with the advantages of the abovementioned zeroing-type models, this paper proposes an adaptive zeroing-type neural dynamics (AZTND) model that embeds adaptive scaling coefficients and adaptive feedback coefficients. For the first time, the AZTND model is applied to the TVQM problem with measuring noise interference.
\par
The remaining part of this paper consists of the following five sections. The problem definition and benchmark scheme are presented in Section \ref{PSF}. The adaptive scale coefficient and adaptive feedback coefficient design framework and the evolution function of the proposed AZTND models are formulated in Section \ref{SMC}. Section \ref{Simulations} contains the corresponding quantitative simulative experiment and results investigation. Finally, the conclusion is arranged in Section \ref{Conclusion}.
Besides, the following parts summarise the main contributions:

\begin{itemize}
\item This paper proposes a novel design framework for constructing the adaptive scale coefficient and adaptive feedback coefficient for the first time, which expedites global convergence and enhances the robustness of the solution system.

\begin{table*}[t]\tiny 
\caption{Comparison between Various Algorithms for TVQM problem (\ref{TVQM})}
\resizebox{14cm}{!}{
\begin{tabular}[l]{@{}l c c c c c c c c c}
	\hline
	&Derivative &Integral  &Adaption &Anti &\multicolumn{4}{c}{MSSRE$^{\ast}$ under Different Noise Conditions}\\
	Model &Information &Information &Control &Perturbations &Noise &Constant &Random &Linear\\
	&Involved &Involved & & &Free &Noise &Noise &Noise\\
	\hline
		Neural network in \cite{IterCompareOne} &\cmark &\xmark &\xmark &\xmark &NA$\dagger$ &NA$\dagger$ &NA$\dagger$ &NA$\dagger$\\
		Adaptive GNN model in \cite{ACGNNCompare}&\cmark &\xmark &\cmark &\xmark &Negligible &Bounded &Bounded &$+\infty$\\
		Original ZNN model in \cite{OZNNCompare} &\cmark &\xmark &\xmark &\xmark &Negligible &Bounded &Bounded &$+\infty$\\
		NCZNN model in \cite{NCZNNCompare} &\cmark &\xmark &\xmark &\cmark &Negligible &Bounded &Bounded &$+\infty$\\
		PTCZNN model in \cite{PTCZNNCompare} &\cmark &\xmark &\xmark &\xmark &Negligible &Bounded &Bounded &$+\infty$\\
		MZNN model in \cite{MZNNCompare} &\cmark &\cmark &\xmark &\cmark &Negligible &Negligible &BS$^\ddagger$ &BS$^\ddagger$\\
		AZTND model (\ref{RACZNN}) &\cmark &\cmark &\cmark &\cmark &Negligible &Negligible &BS$^\ddagger$ &BS$^\ddagger$\\\hline
\end{tabular}\label{RiccatiCompare}}

	\noindent{\footnotesize{$\ast$ Note that MSSRE denotes the maximal steady-state residual errors.}}\\
	\noindent{\footnotesize{$^\dagger$NA indicates that the item does not apply to the algorithm or model in the corresponding papers.}}\\
	\noindent{\footnotesize{$^\ddagger$BS indicates that the maximal steady-state residual error of the corresponding situation is bounded tightly.}}
\end{table*}

\item Based on the design framework, an AZTND model for solving the TVQM problem with perturbed measurement noise is proposed. Subsequently, the global convergence of the AZTND model has analyzed the Lyapunov stability theory.
\item Corresponding quantitative numerical experiments are given to substantiate the performance of the AZTND model applied to the TVQM problem with various measurement noise pollution.
\item A dynamic localization scheme is proposed based on the AZTND model with adaptive coefficients, which has superior robustness and solution accuracy to existing schemes.
\end{itemize}

\section{Problem Definition and Related Scheme}\label{PSF}
The typical form of the TVQM problem is presented as
\begin{equation}\label{TVQM}
	\text{min}~\frac{1}{2}\vec{z}^{\text{T}}(t)M(t)\vec{z}(t)+\vec{b}^{\text{T}}(t)\vec{z}(t),
\end{equation}
where parameters $M(t)\in\mathbb{R}^{n\times n}$ and $\vec{b}(t)\in\mathbb{R}^{n}$ mean the smoothly time-varying Hessian matrix and vector, respectively. The parameter $\vec{z}(t)\in\mathbb{R}^{n}$ represents the unknown vector that should be solved online. The transpose symbol is the superscript $^{\text{T}}$. For further investigation and solving the TVQM problem (\ref{TVQM}), a function $F(\vec{z}(t),t)$ is defined as $F(\vec{z}(t),t)=\frac{1}{2}\vec{z}^{\text{T}}(t)M(t)\vec{z}(t)+\vec{b}^{\text{T}}(t)\vec{z}(t)$. Consequently, the gradient of the function $F(\vec{z}(t),t)$ is described as
\begin{equation}\label{FGrad}
	\nabla F(\vec{z}(t),t)=\frac{\partial F(\vec{z}(t),t)}{\partial \vec{z}(t)}=M(t)\vec{z}(t)+\vec{b}(t).
\end{equation}
Noteworthily, by zeroing $\nabla F(\vec{z}(t),t)$ in each time instant $t\in [0,+\infty]$, the theoretical solution of TVQM problem (\ref{TVQM}) can be obtained in real-time. Hence, the following equation is formulated as $M(t)\vec{z}(t)+\vec{b}(t)=0$.
The afterward error function is arranged to tune the evolution direction of the solving system:
\begin{equation}\label{ErrFun}
	\vec{\epsilon} (t)=M(t)\vec{z}(t)+\vec{b}(t).
\end{equation}
According to the OZNN model construction framework, the evolution direction of the error function (\ref{ErrFun}) should be satisfied that $\vec{\dot\epsilon}(t)=-\eta\Omega (\vec{\epsilon}(t))$, where $\eta$ represents the scale coefficient and $\Omega(\cdot): \mathbb{R}^{n} \to \mathbb{R}^{n}$ denotes the activation function. Therefore, the OZNN model employed for the TVQM problem (\ref{TVQM}) is described as
\begin{equation}
	M(t)\vec{\dot z}(t)=-\dot M(t)\vec{z}(t)-\vec{\dot b}(t)-\eta\Omega\big{(}M(t)\vec{z}(t)+\vec{b}(t)\big{)},
\end{equation}
where parameters $\dot M(t)$, $\vec{\dot{z}}(t)$, and $\vec{\dot{b}}(t)$ represent time derivatives of $M(t)$, $\vec{z}(t)$, and $\vec{b}(t)$, respectively.
\par
The performance comparison among the existing algorithms and the proposed AZTND model when solving the TVQM problem (\ref{TVQM}) are arranged in Table \ref{RiccatiCompare}.

\section{AZTND Model Construction}\label{SMC}
Various perturbed measurement noises downgrade the accuracy of the solution system and even lead to collapse. Thus, this paper proposes an AZTND model to enhance the robustness under various noise interference. Besides, to avoid the inflexibility of manually setting the scale coefficient and to speed up the convergence of the solving system, this paper formulates a more flexible scale coefficient construction method termed residual-based adaptive scale coefficient.
\par
The evolution direction of the error function $\vec{\epsilon}(t)$ in the AZTND model is formulated as
\begin{equation}
	\vec{\dot \epsilon}(t)=-\xi(\vec{\epsilon}(t))\vec{\epsilon}(t)-\kappa(\epsilon(t))\int_0^t\vec{\epsilon}(\delta)\text{d}\delta,
\end{equation}
where parameters $\xi(\cdot)> 0: \mathbb{R}^{n}\to \mathbb{R}$ and $\kappa(\cdot)> 0: \mathbb{R}^{n}\to \mathbb{R}$ represent the adaptive scale and feedback coefficient, respectively. The following method can be employed to construct the adaptive scale coefficient $\xi(\cdot)$:
\begin{itemize}
	\item Power adaptive scale coefficient:
	\begin{equation}\label{ASCDef}
		\xi(\epsilon(t)) = \|\vec{\epsilon}(t)\|_{\text{2}}^{\eta} + a,
	\end{equation}
	where the parameter $a > 1$.
\end{itemize}
Further, the following method can be adopted to implement the adaptive feedback coefficient $\kappa(\cdot)$:
\begin{itemize}
		\item Power adaptive feedback coefficient:
	\begin{equation}\label{AFCDef}
		\kappa(\vec{\epsilon}(t)) = \|\int_0^t\vec{\epsilon}(\delta)\text{d}\delta\|^b_{\text{2}}+c,
	\end{equation}
	where the parameter $b>0$ and $c>0$.
\end{itemize}
Consequently, the proposed AZTND model with adaptive scale coefficient for solving the TVQM problem (\ref{TVQM}) is written as follows:
\begin{eqnarray}\label{RACZNN}
	\begin{split}
		 M(t)\vec{\dot z}(t)=&-\dot M(t)\vec{z}(t)-\vec{\dot b}(t)\\
		&-\xi(\vec{\epsilon}(t))\big{(}M(t)\vec{z}(t)+\vec{b}(t)\big{)} \\
		&-\kappa(\vec{\epsilon}(t))\int_{0}^{t}(M(\delta)\vec{z}(\delta)+\vec{b}(\delta))\text{d}\delta.
	\end{split}
\end{eqnarray}
Besides, the AZTND model (\ref{RACZNN}) inevitably is perturbed by various measurement noises in the practical application. Thereupon, the AZTND model (\ref{RACZNN}) for solving TVQM problem (\ref{TVQM}) perturbed by noise is described as
\begin{eqnarray}\label{RACZNNNoise}
	\begin{split}
		 M(t)\vec{\dot z}(t)=&-\dot M(t)\vec{z}(t)-\vec{\dot b}(t) \\
		&-\xi(\vec{\epsilon}(t))\big{(}M(t)\vec{z}(t)+\vec{b}(t)\big{)}\\
		&-\kappa(\vec{\epsilon}(t))\int_{0}^{t}(M(\delta)\vec{z}(\delta)+\vec{b}(\delta))\text{d}\delta + \vec{\vartheta}(t),
	\end{split}
\end{eqnarray}
where the noise perturbation item $\vec{\vartheta}(t)\in \mathbb{R}^{n}$.
\par
Taking into account that convergence is a key criterion for the AZTND model (\ref{RACZNN}), we propose the following theorem and corresponding proof process to analyze the global convergence of the AZTND model (\ref{RACZNN}).
\par
{\it Theorem 1:} For any solvable TVQM problem (\ref{TVQM}), the proposed AZTND model (\ref{RACZNN}) globally converges to zero from any random initial state.
\par

{\it Proof:} The $i$th subsystem of the AZTND model evolution function $\vec{\dot \epsilon}(t)=-\xi(\vec{\epsilon}(t))\vec{\epsilon}(t)-\kappa(\vec\epsilon(t))\int_0^t\vec{\epsilon}(\delta)\text{d}\delta$ is written as
\begin{equation}\label{SubErr}
	\dot\epsilon_{i}(t)=-\xi\big{(}\epsilon(t)\big{)}\epsilon_{i}(t)-\kappa(\epsilon(t))\int_{0}^{t}\epsilon_{i}(\delta)\text{d}\delta,
\end{equation}
The following Lyapunov function is presented for investigating the global convergence of the system (\ref{SubErr}):
\begin{equation}\label{Lya}
	y_{i}(t)=\epsilon_{i}^2(t)+\kappa(\epsilon(t))\big{(}\int_{0}^{t}\epsilon_{i}(\delta)\text{d}\delta\big{)}^2/2\ge 0,
\end{equation}
which indicates the function $y_{i}(t)>0$ when $\epsilon_{i}(t)\neq 0$ or $\int_{0}^{t}\epsilon_{i}(\delta)\text{d}\delta \neq 0$. If and only if $\epsilon_{i}(t)=\int_{0}^{t}\epsilon_{i}(\delta)\text{d}\delta=0$, $y_{i}(t)=0$. Thus, the Lyapunov function $y_{i}(t)$ is positive semi-definite. Considering that the adaptive feedback coefficient $\kappa(\epsilon(t))$ is a constant $\kappa$ in each time interval and taking the time derivative of the function (\ref{Lya}) leads:

\begin{eqnarray*}
	\begin{split}
		\frac{\text{d}y_{i}(t)}{\text{d}t}&=\epsilon_{i}(t)\dot\epsilon_{i}(t)+\kappa(\vec{\epsilon}(t))\epsilon_{i}(t)\int_{0}^{t}\epsilon_{i}(\delta)\text{d}\delta \\
		&=\epsilon_{i}(t)(\dot\epsilon_{i}(t)+\kappa(\vec{\epsilon}(t))\int_{0}^{t}\epsilon_{i}(\delta)\text{d}\delta)\\
		&=-\xi(\vec\epsilon(t))\epsilon_{i}^2(t)\leq 0.
	\end{split}
\end{eqnarray*}
That is to say, the Lyapunov function $\dot y_{i}(t)$ is negative semi-definite. Thus, according to the definition of the Lyapunov theory, the function $\epsilon_{i}(t)$ will ultimately converge to zero. It can be generalized and concluded that $\epsilon_{i}(t)$ globally converges to zero for each $i\in 1,2,...,n$. In summary, the error function $\vec\epsilon(t)$ global converges to zero over time. In other words, the proposed AZTND model (\ref{RACZNN}) globally converges to the theoretical solution of the TVQM problem (\ref{TVQM}).
\par
The proof is thus completed. $\hfill\blacksquare$

\begin{figure}[htbp]\centering
\subfigure[]{\includegraphics[scale=0.4]{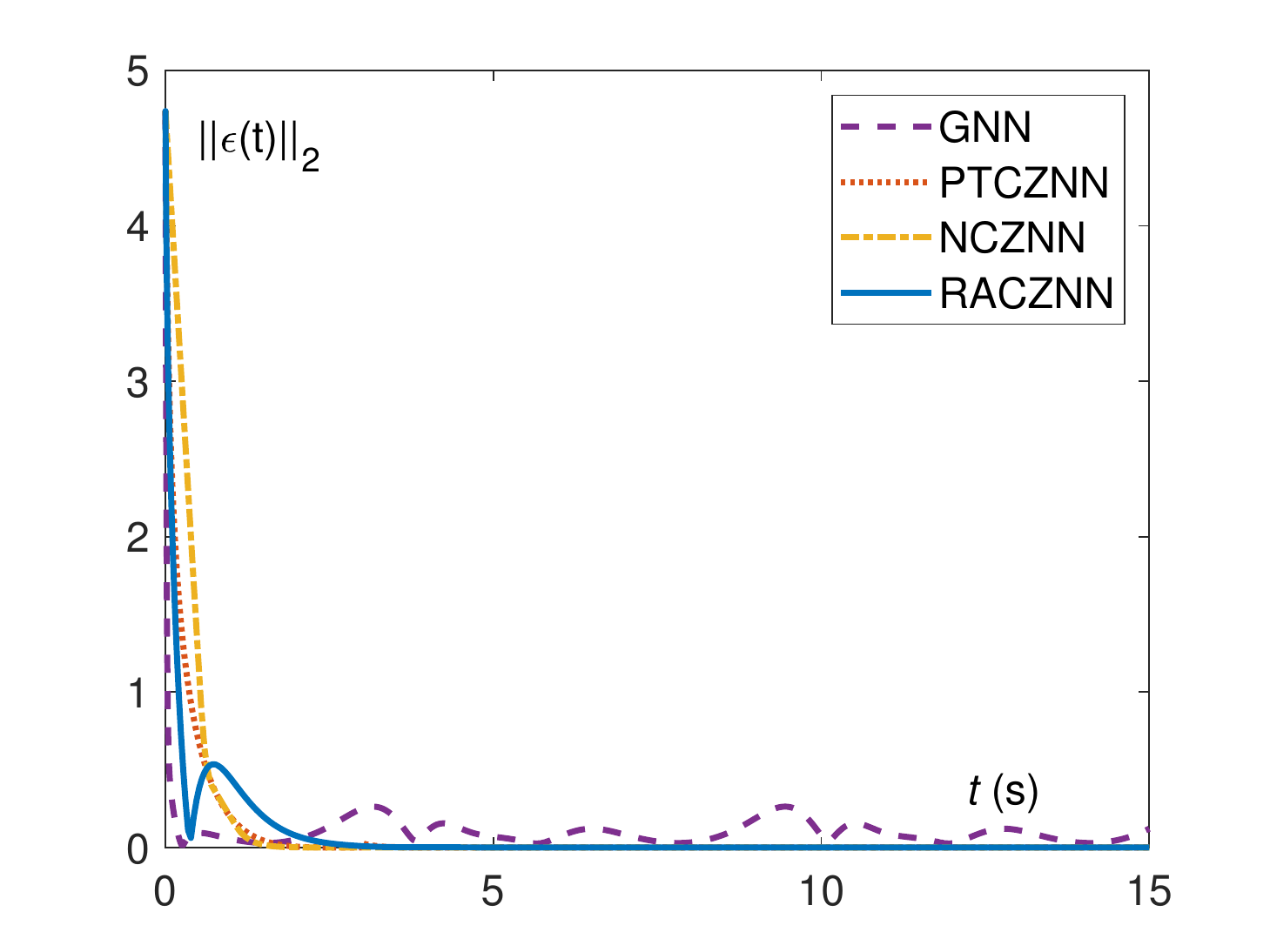}}
\subfigure[]{\includegraphics[scale=0.4]{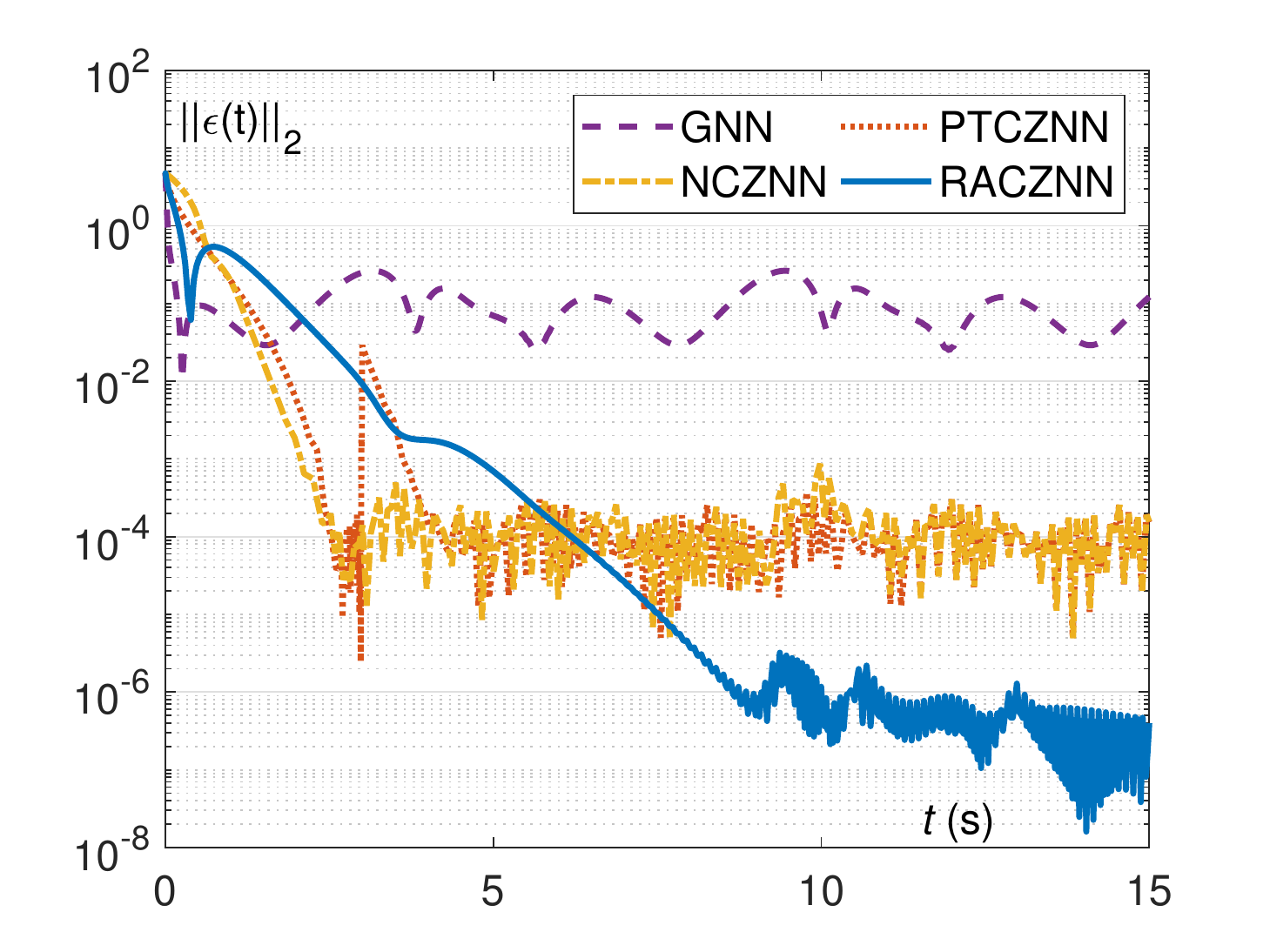}}
\caption{Performance of the GNN (\ref{GNNCompare}), PTCZNN (\ref{PTCZNNCompare}), NCZNN (\ref{NCZNNCompare}), and AZTND (\ref{RACZNN}) are applied to noise-free TVQM problem (\ref{TVQM}) of the Example (\ref{EA}). (a) Denoting the residual error $||\vec{\epsilon}(t)||_{\text{2}}$ of models. (b) The logarithm of residual error $||\vec{\epsilon}(t)||_{\text{2}}$.
}
\label{FreeNorm}
\end{figure}

\section{Simulations}\label{Simulations}
Experiments are designed and performed in this section. First, the simulation of the AZTND model (\ref{RACZNN}) applied to the TVQM problem (\ref{TVQM}) is concluded and visualized. Secondly, we compare the performance of the AZTND model (\ref{RACZNN}) with other state-of-the-art neural network models, specifically the gradient-based neural network (GNN) model, predefined-time convergent ZNN (PTCZNN) model, nonconvex and bound constraint ZNN (NCZNN) model.
\subsection{Example 1: Time-varying Situation}\label{EA}
In this simulation, the time-varying matrix and vector in the TVQM problem (\ref{TVQM}) are constructed as follows:
\begin{eqnarray*}
    M(t)=
    \begin{bmatrix}
        0.5\text{sin}(t)+2& \text{cos}(t)\\
        \text{cos}(t) & 0.5\text{cos}(t)+2
    \end{bmatrix},
    \vec{b}(t)=
    \begin{bmatrix}
        \text{sin}(t)\\
        \text{cos}(t)
    \end{bmatrix}.
\end{eqnarray*}

\begin{figure}[htbp]\centering
\subfigure[]{\includegraphics[scale=0.4]{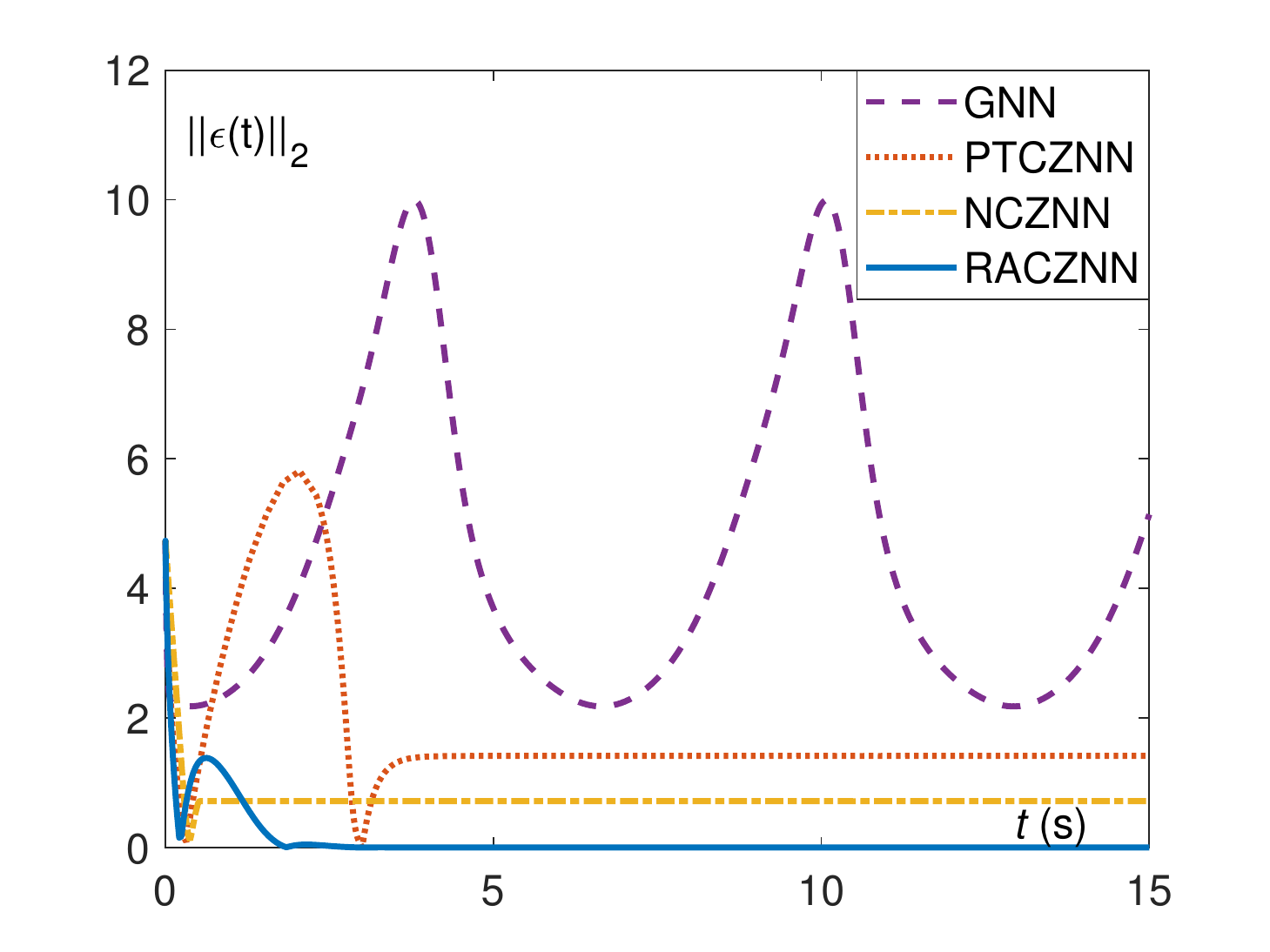}}
\subfigure[]{\includegraphics[scale=0.4]{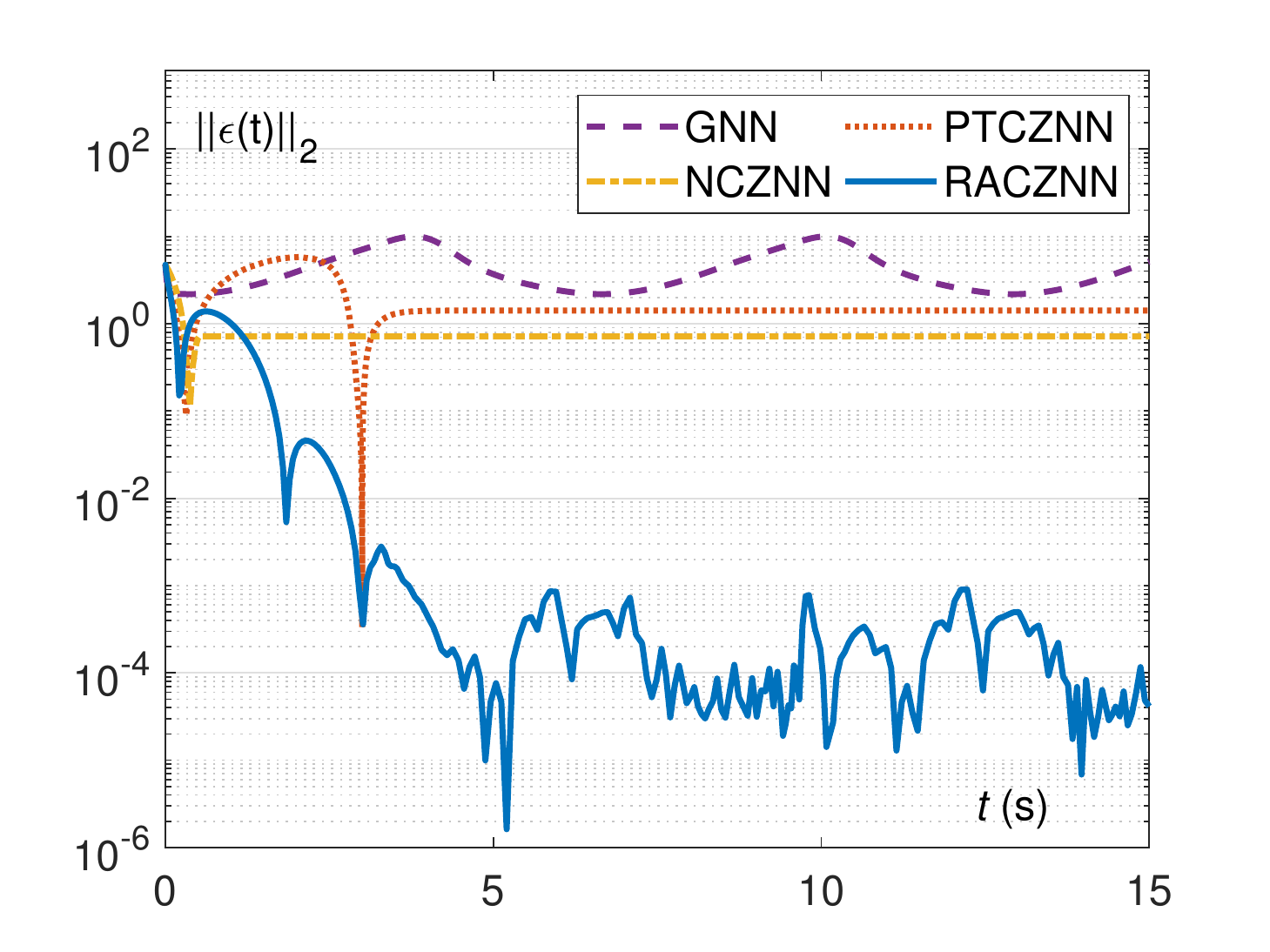}}
\caption{Performance comparison among different models for solving the TVQM problem (\ref{TVQM}) of the Example (\ref{EA}) with constant noise $\vec{\vartheta}(t)=\vec{\vartheta}=[5]^2$. (a) Residual error $||\vec{\epsilon}(t)||_{\text{2}}$ of models. (b) The logarithm of the residual error $||\vec{\epsilon}(t)||_{\text{2}}$.
}
\label{ConNorm}
\end{figure}
\begin{figure}[htbp]\centering
\subfigure[]{\includegraphics[scale=0.4]{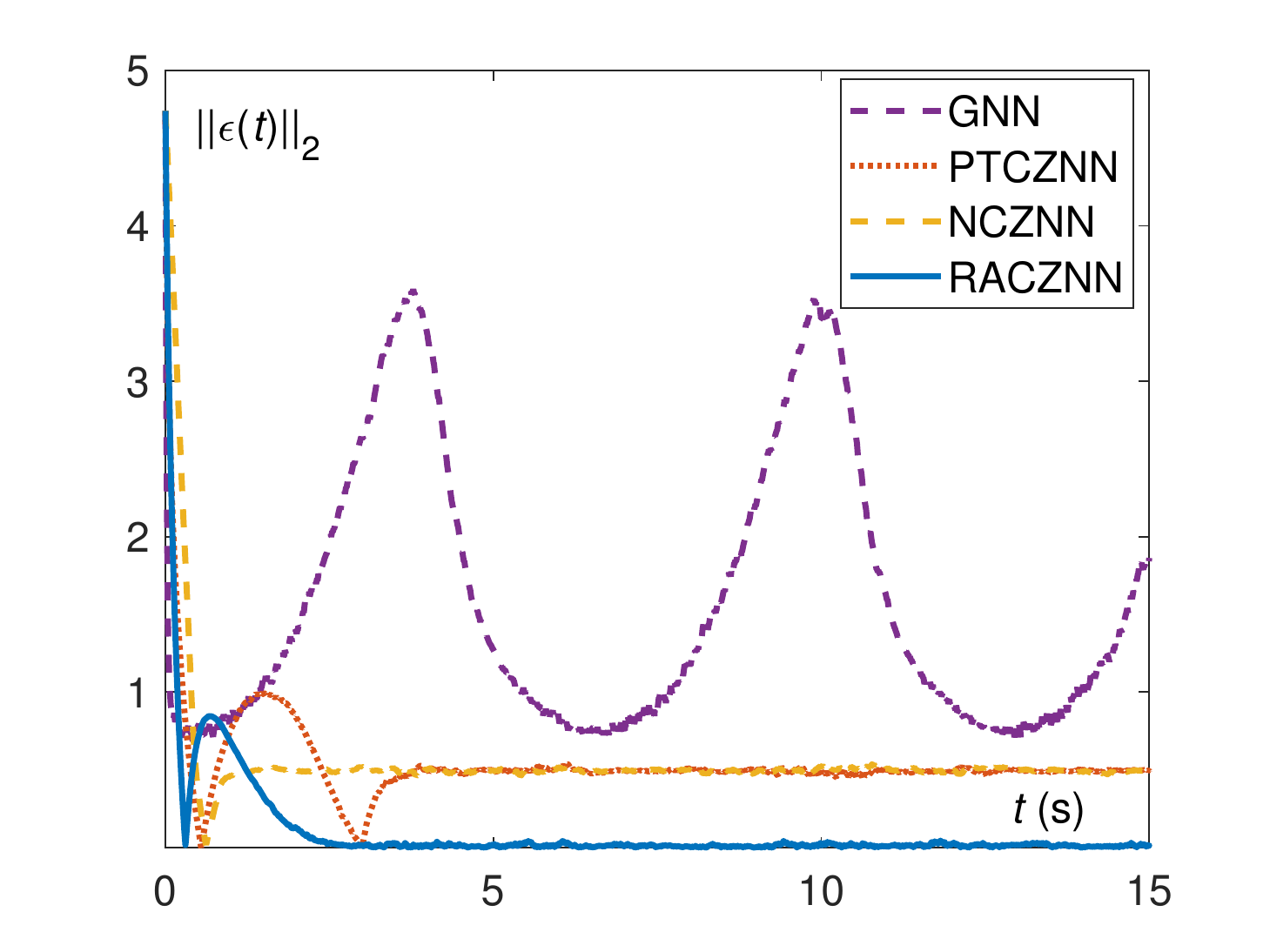}}
\subfigure[]{\includegraphics[scale=0.4]{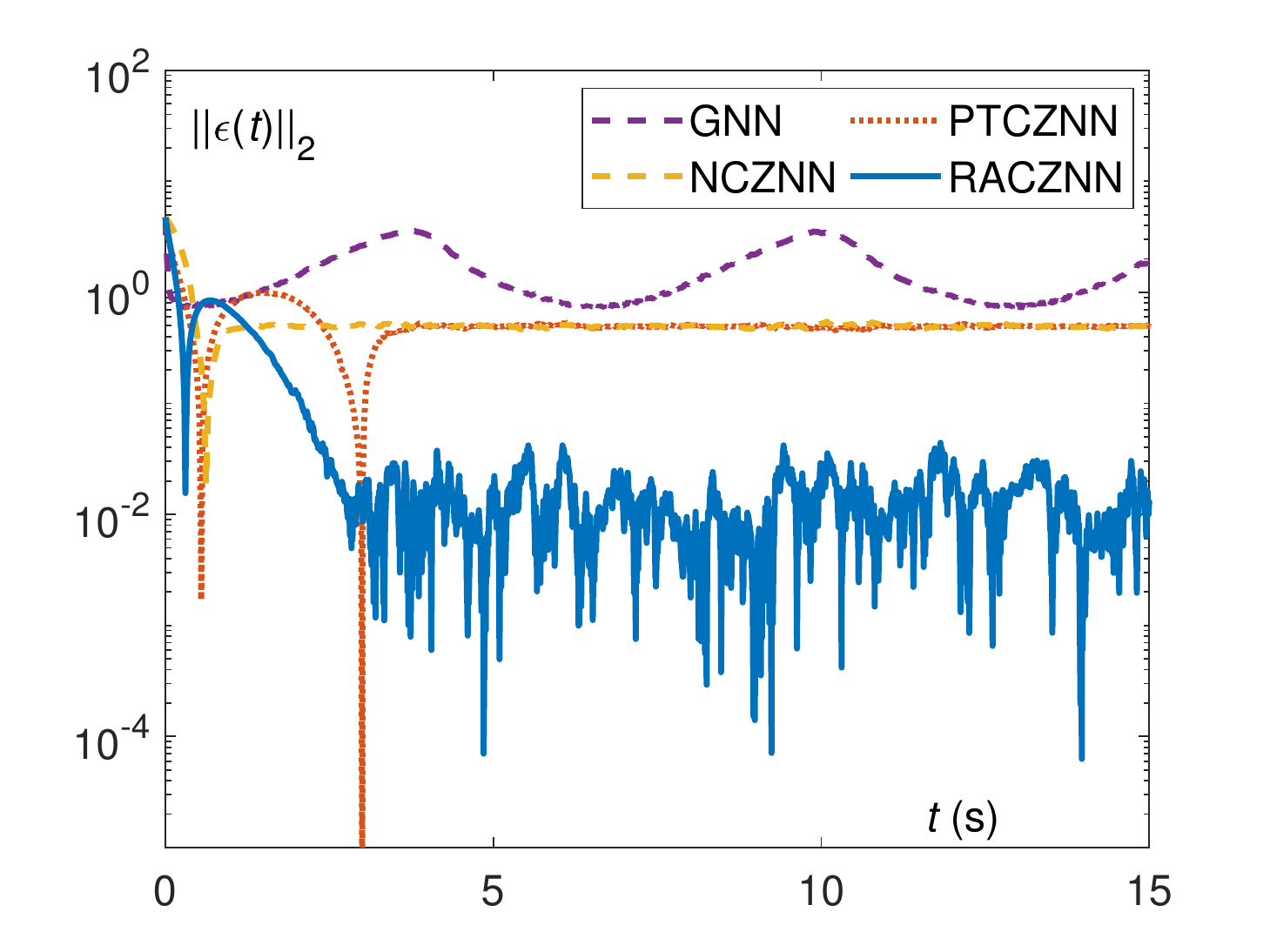}}
\caption{The performance of different models with bounded random noise $\vec{\vartheta}(t)=\vec{\varrho}(t)\in [0.5, 3]^2$. (a) Denoting the residual error $||\vec{\epsilon}(t)||_{\text{2}}$ of the models. (b) Representing the logarithm of residual error $||\vec{\epsilon}(t)||_{\text{2}}$.
}\label{RandNorm}
\end{figure}

The adaptive scale coefficient $\xi(\cdot)$ and adaptive feedback coefficient $\kappa(\cdot)$ of the proposed AZTND model (\ref{RACZNN}) are set as $\xi(\vec{\epsilon}(t)) = ||\vec{\epsilon}(t)||_{\text{2}}^{3} + 5$ and $\kappa(\vec{\epsilon}(t)) = 5^{||\int_0^t\vec{\epsilon}(\delta)\text{d}\delta||_{\text{2}}}+5$, respectively. The corresponding quantitative simulation results of the example (\ref{EA}) are arranged in Figures \ref{FreeNorm} to \ref{RandNorm}. Besides, the following models are introduced to solve the TVQM problem (\ref{TVQM}) as a comparison of the AZTND model (\ref{RACZNN}).

\begin{itemize}
\item The GNN model is presented in \cite{GNNCompare}.
\begin{eqnarray}\label{GNNCompare}
	\begin{split}
	\vec{\dot z}(t)=- \gamma M^{\text{T}}(t)(M(t)\vec{z}(t)+\vec{b}(t)).
	\end{split}
\end{eqnarray}
\begin{figure}[htbp]\centering
\subfigure[]{\includegraphics[scale=0.4]{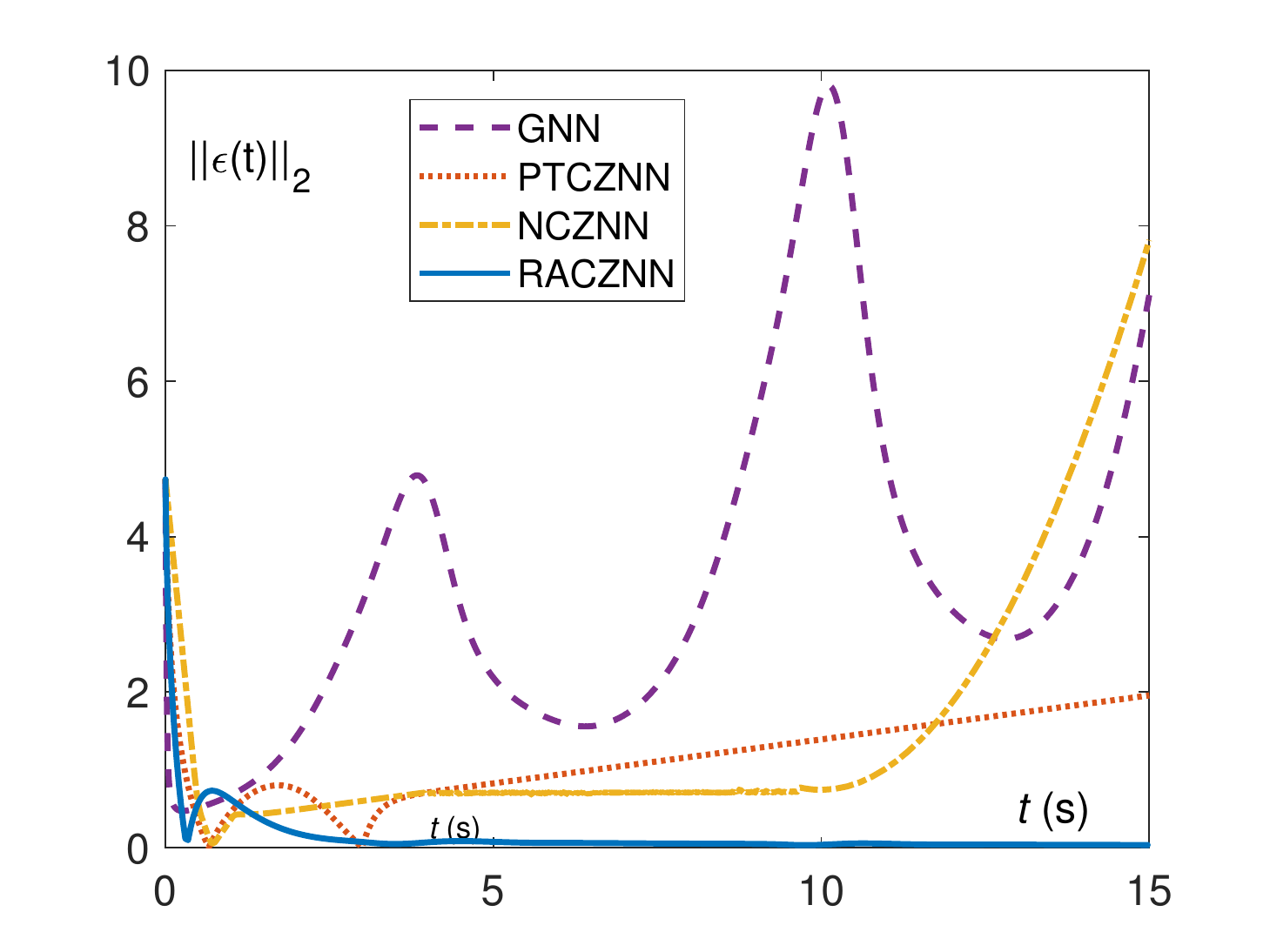}}
\subfigure[]{\includegraphics[scale=0.4]{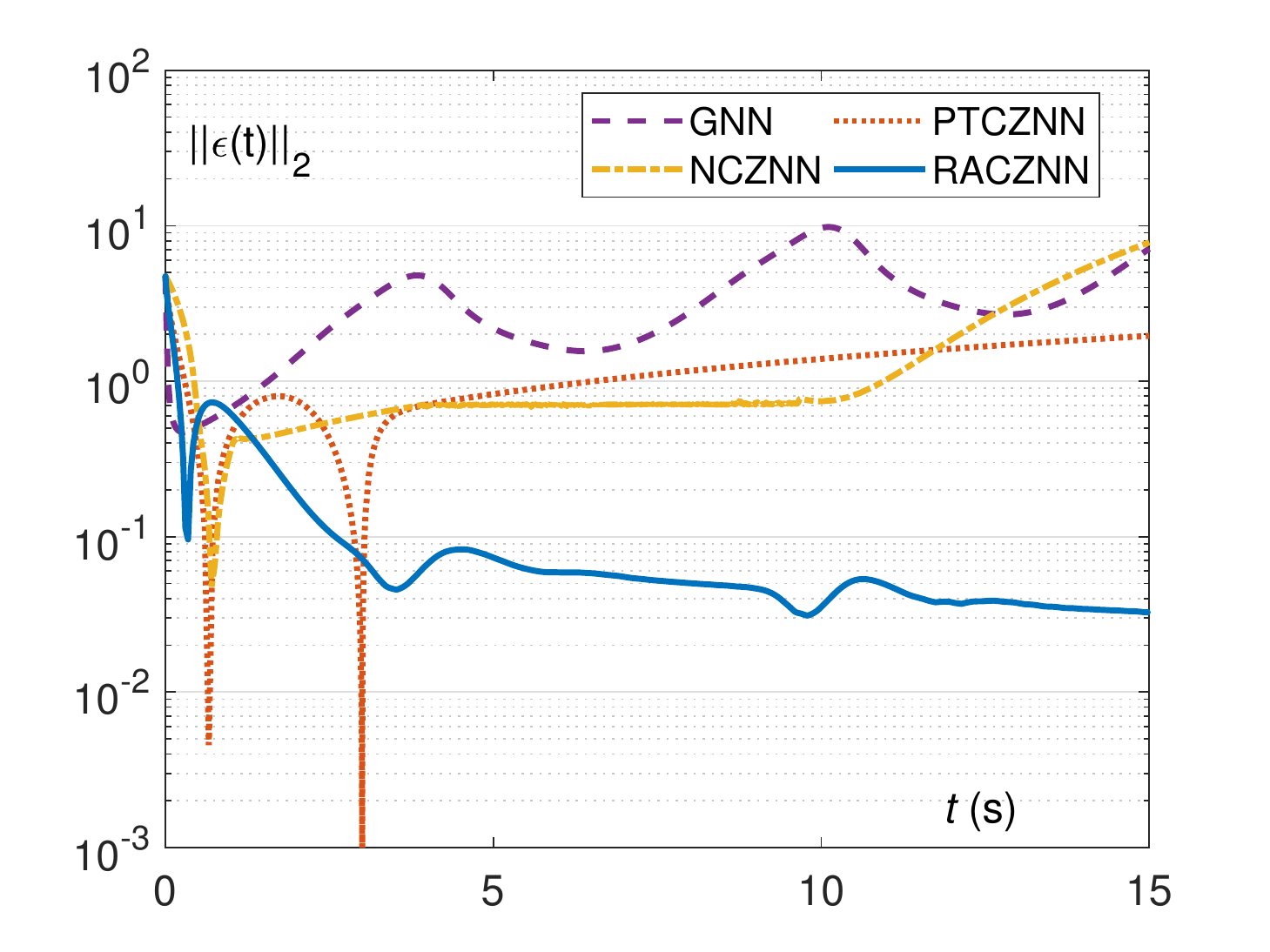}}
\caption{
The performance of different models under linear noise $\vec{\vartheta}(t)\in \mathbb{R}^{n^2}$ with each subelement be set as $0.4\times t$. (a) Denoting the residual error $||\vec{\epsilon}(t)||_{\text{2}}$ of the models. (b) Representing the logarithm of residual error $||\vec{\epsilon}(t)||_{\text{2}}$.
}
\label{LinearNorm}
\end{figure}
\item PTCZNN model is presented in \cite{PTCZNNCompare}.
\begin{eqnarray}\label{PTCZNNCompare}
	\begin{split}
	M(t)\vec{\dot z}(t)=&-\dot M(t)\vec{z}(t)-\vec{\dot b}(t)\\
	&-\gamma \frac{\text{exp}(t)-1}{(t_c-t)\text{exp}(t)}(M(t)\vec{z}(t)-\vec{d}(t)).
	\end{split}
\end{eqnarray}

\item NCZNN model is presented in \cite{NCZNNCompare}.
\begin{eqnarray}\label{NCZNNCompare}
	\begin{split}
	M(t)\vec{\dot z}(t)=&-\dot M(t)\vec{z}(t)-\vec{\dot b}(t)\\
	&-\gamma R_{\Upsilon}(M(t)\vec{z}(t)-\vec{d}(t)),
	\end{split}
\end{eqnarray}
where the parameter $R_{\Upsilon}(\cdot)$ represents the non-convex and bounded activation function.
\end{itemize}
Note that the scale parameter $\gamma$ of the GNN model (\ref{GNNCompare}), PTCZNN model (\ref{PTCZNNCompare}), and NCZNN model (\ref{NCZNNCompare}) are arranged as 5.

\begin{figure*}[htbp]\centering
\subfigure[]{\includegraphics[scale=0.35]{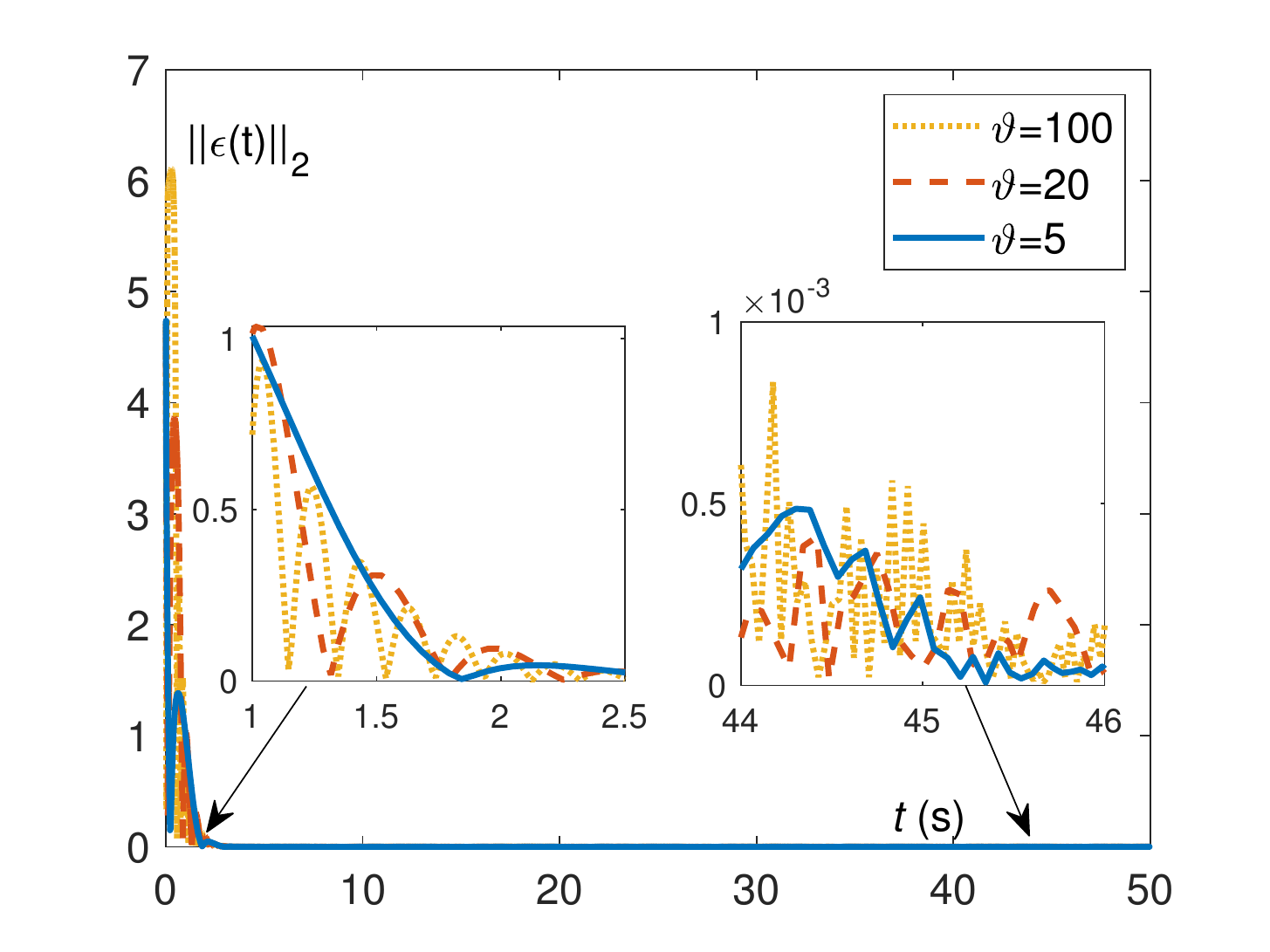}}
\subfigure[]{\includegraphics[scale=0.35]{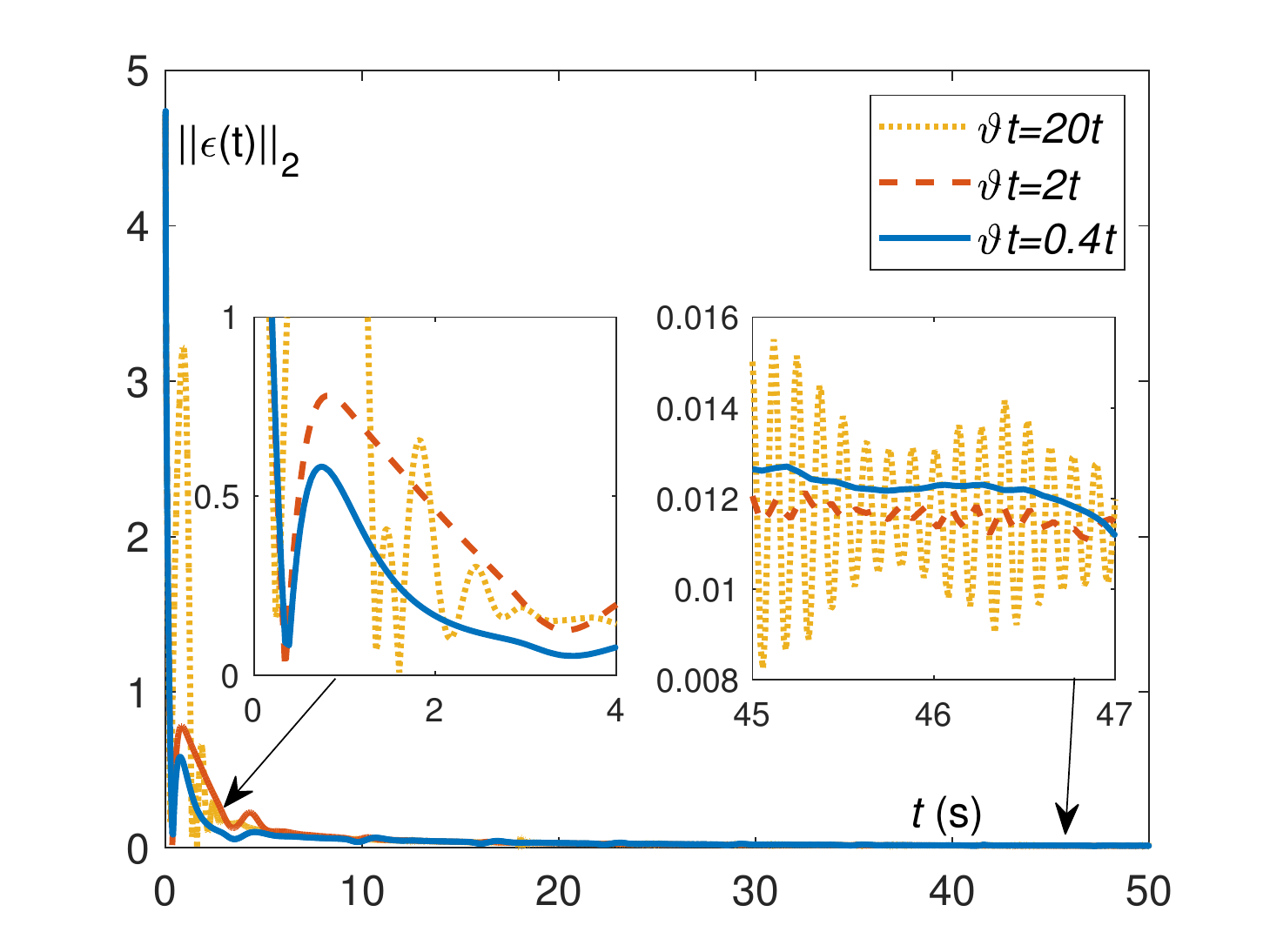}}
\subfigure[]{\includegraphics[scale=0.35]{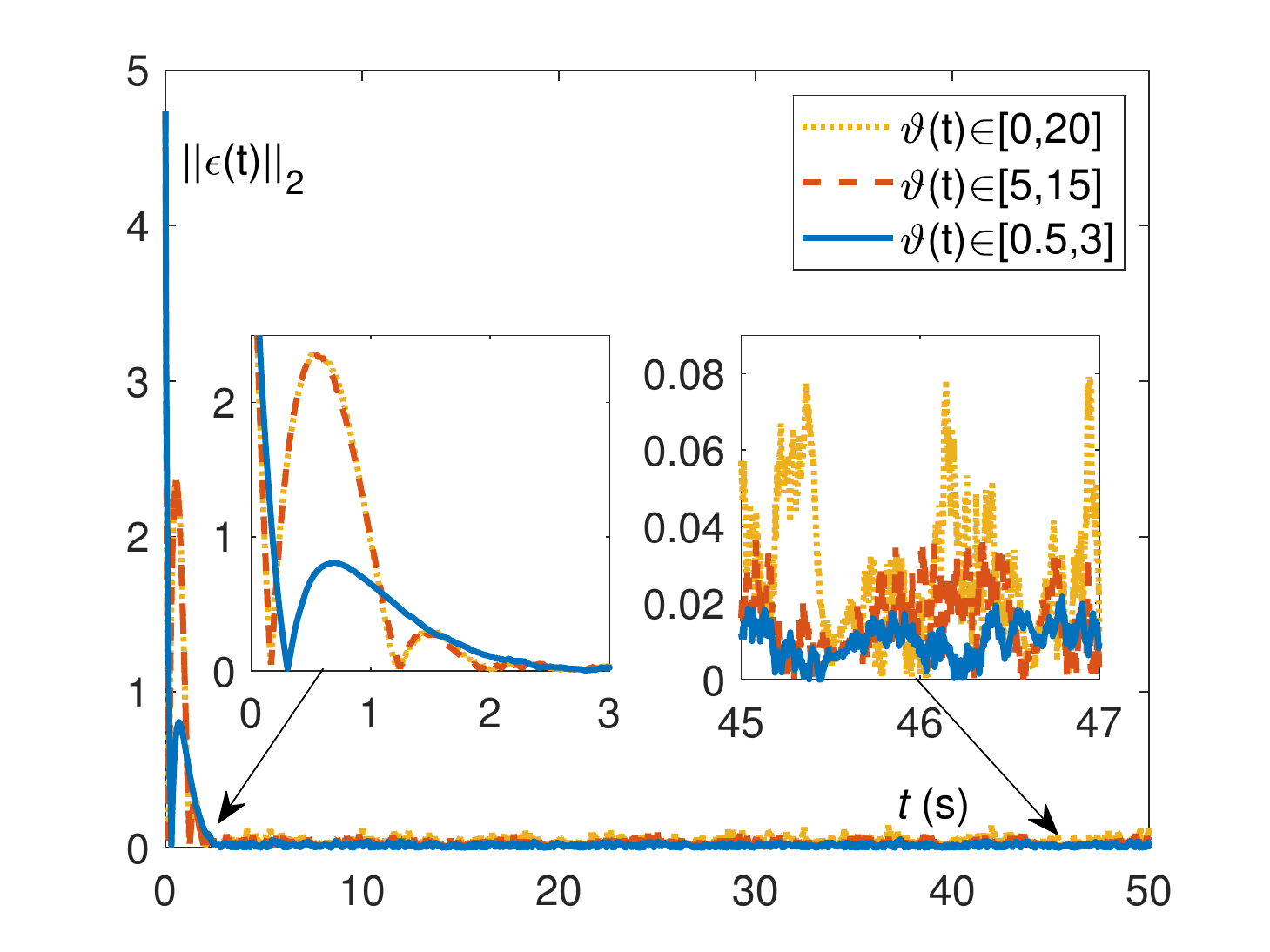}}
\subfigure[]{\includegraphics[scale=0.35]{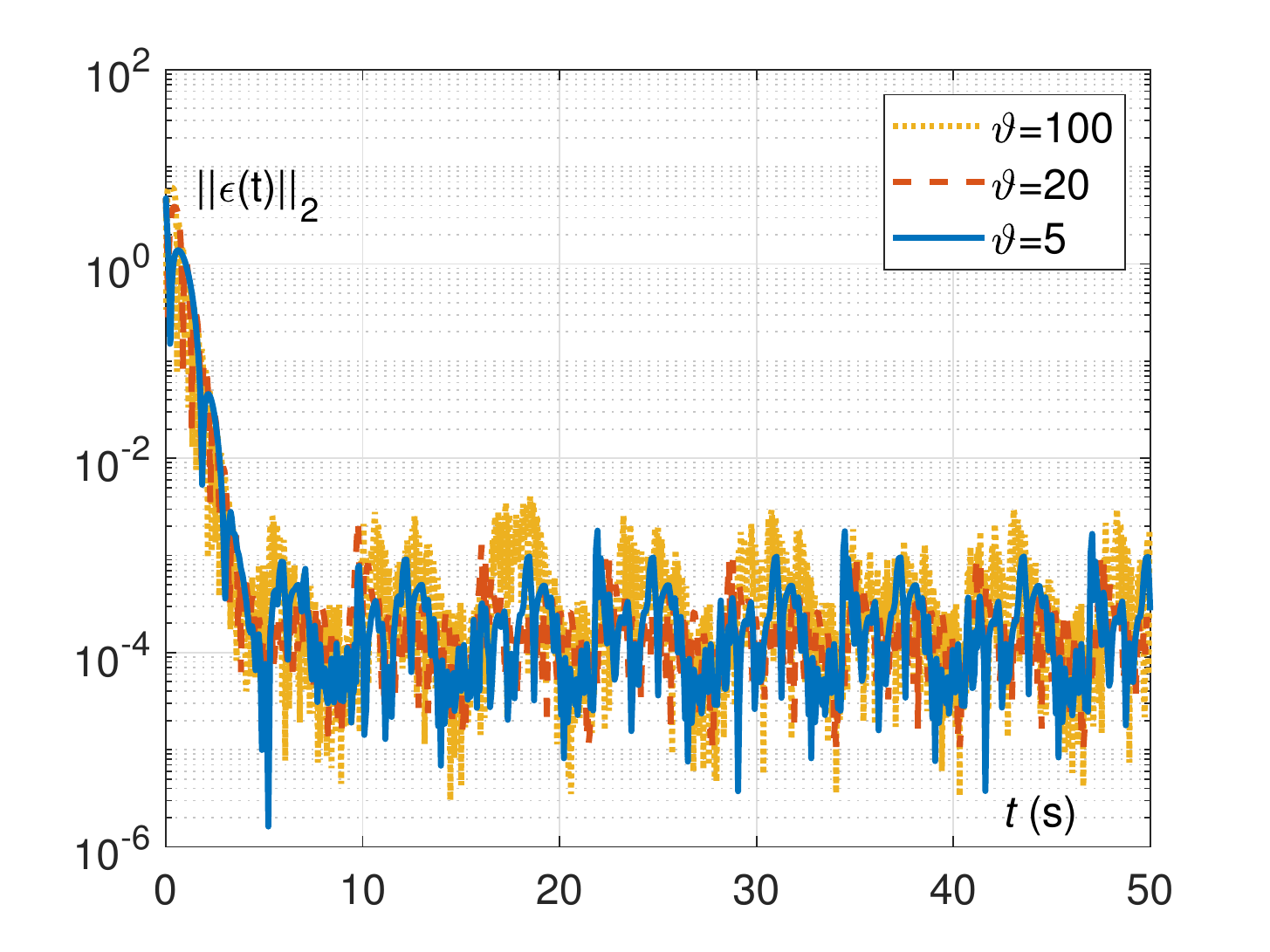}}
\subfigure[]{\includegraphics[scale=0.35]{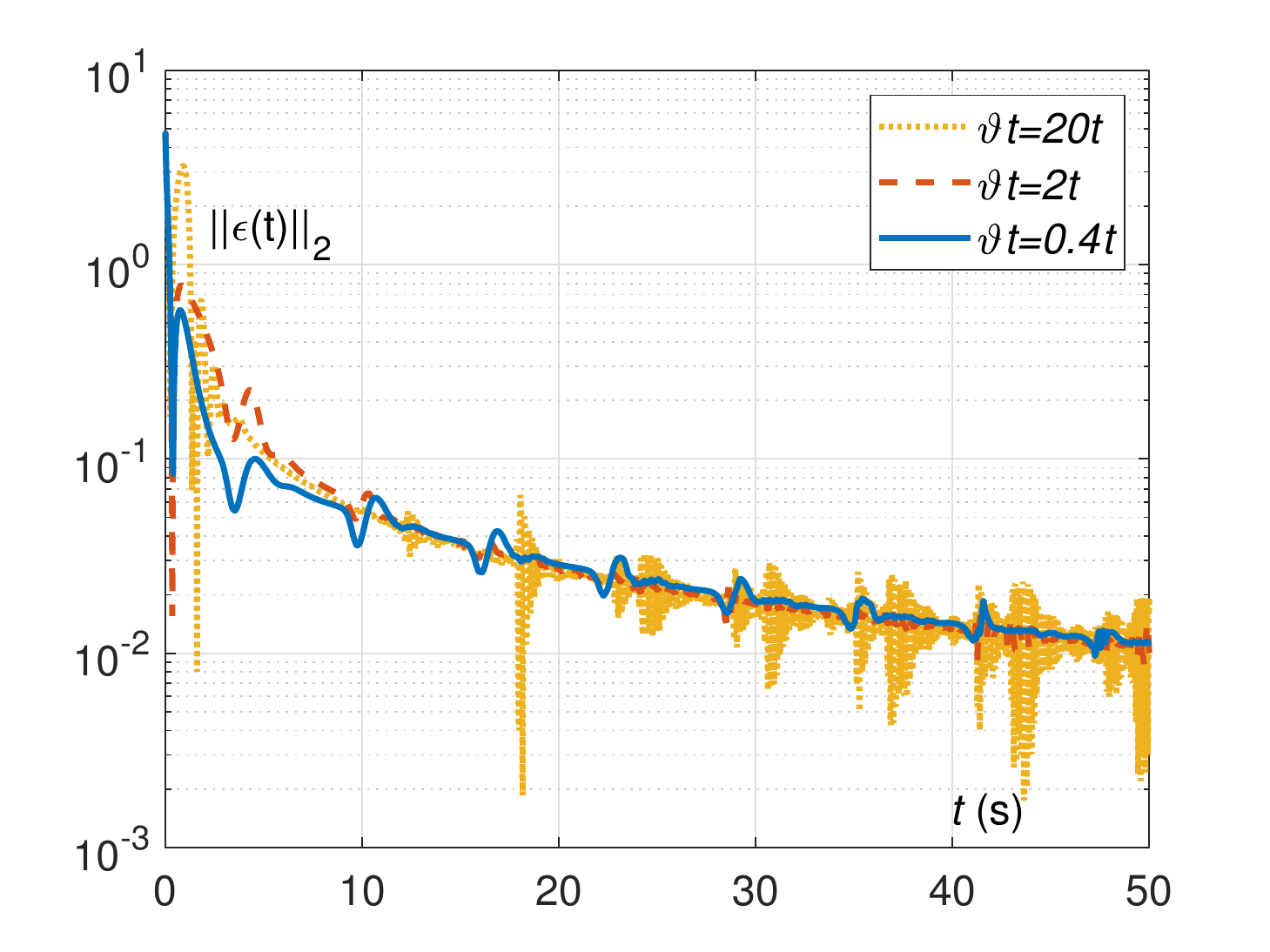}}
\subfigure[]{\includegraphics[scale=0.35]{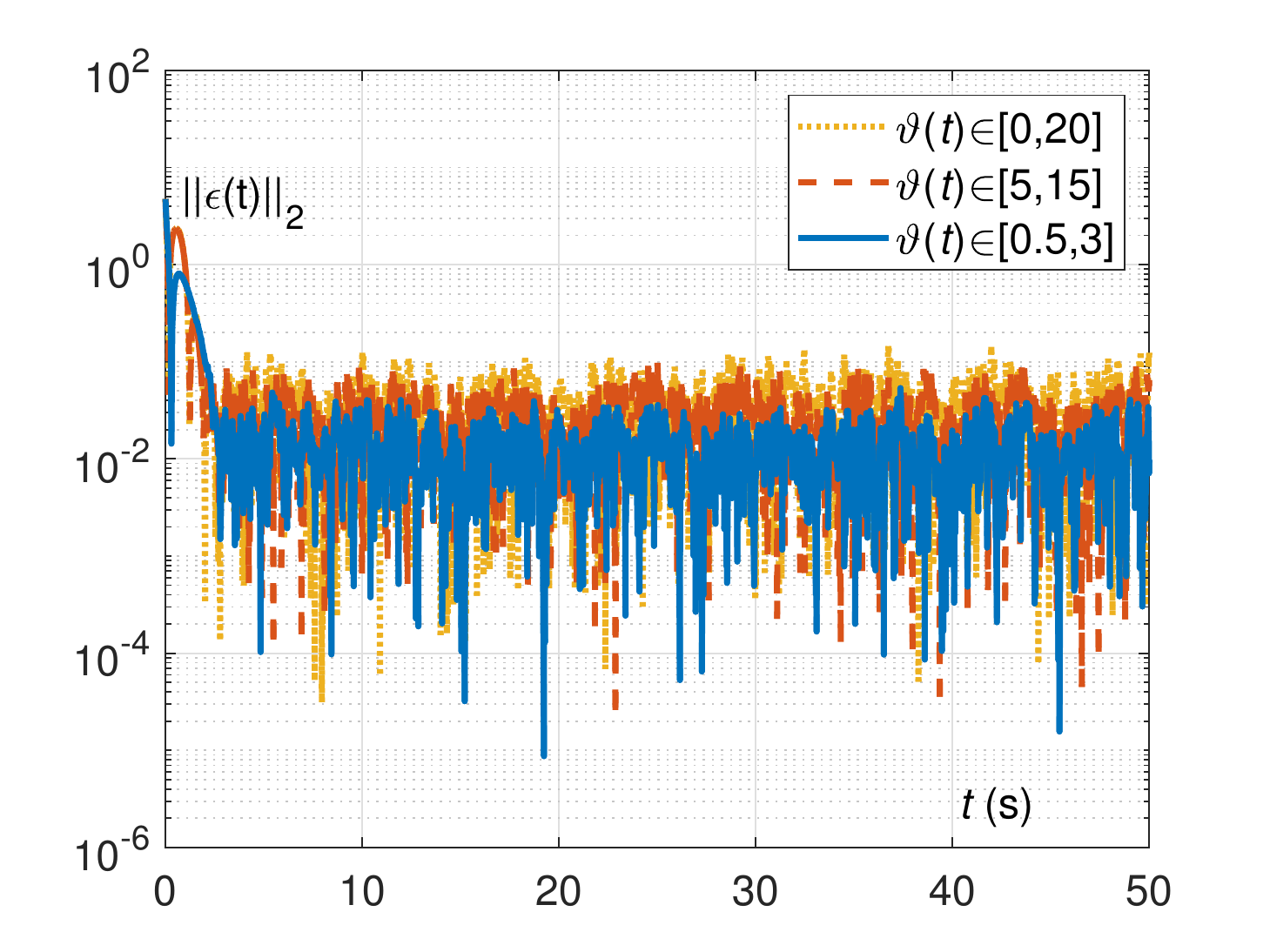}}
\caption{Robustness performance of the proposed AZTND model (\ref{RACZNN}) under noise interference cases. (a) and (d) representing the residual error $||\vec{\epsilon}(t)||_{\text{2}}$ of the AZTND model (\ref{RACZNN}) perturbed by constant noise $\vartheta=5$, $\vartheta=20$, and $\vartheta=100$, respectively. (b) and (e) denoting the residual error $||\vec{\epsilon}(t)||_{\text{2}}$ of the AZTND model (\ref{RACZNN}) in case of time-varying linear noise $\vartheta (t)=0.4t$, $\vartheta (t)=2t$, and $\vartheta(t)=20t$, respectively. (c) and (f) denoting the residual error $||\vec{\epsilon}(t)||_{\text{2}}$ of the AZTND model (\ref{RACZNN}) perturbed by bounded random noise $\vartheta(t)=\varrho(t)\in[0.5,3]$, $\vartheta(t)=\varrho(t)\in[0.5,3]$, and $\vartheta(t)=\varrho(t)\in[0.5,3]$, respectively.
}
\label{Compare}
\end{figure*}

\subsubsection{AZTND Model without noise}
The quantitative experiment visualization results synthesized by the AZTND model (\ref{RACZNN}) for solving the TVQM problem example (\ref{EA}) in the noise-free case are arranged in Figure \ref{FreeNorm}. As demonstrated in Figure \ref{FreeNorm} (a), beginning with a randomly generated initial vector-formed value, the residual error $||\vec{\epsilon}(t)||_{\text{2}}$ of the proposed AZTND model (\ref{RACZNN}) sharply approaches zero, which means the solving system globally converges to the theoretical solution. Among the compared five models, the AZTND model (\ref{RACZNN}) has the second convergence speed. Logarithms of the models' residual error $||\vec{\epsilon}(t)||_{\text{2}}$ are shown in Figure \ref{FreeNorm} (b) which depicts the models' solution accuracy. As observed in Figure \ref{FreeNorm} (b), the proposed AZTND model (\ref{RACZNN}) has significantly higher accuracy when solving the noise-free TVQM problem (\ref{EA}) compared with the GNN model (\ref{GNNCompare}), PTCZNN model (\ref{PTCZNNCompare}), and NCZNN model (\ref{NCZNNCompare}). Specifically, the GNN model (\ref{GNNCompare}) converges to order $10^{-1}$, PTCZNN model (\ref{PTCZNNCompare}) and NCZNN model (\ref{NCZNNCompare}) converge to order $10^{-3}$, and the proposed AZTND model (\ref{RACZNN}) converges to order $10^{-5}$. Furthermore, the AZTND model (\ref{RACZNN}) converges to the steady-state residual error at $8.5$ s. Compared with state-of-the-art models, the proposed AZTND model (\ref{RACZNN}) has a competitive robustness and convergence speed performance.

\begin{figure*}[htbp]\centering
\subfigure[Theoretical and estimated trajectory]{\includegraphics[scale=0.33]{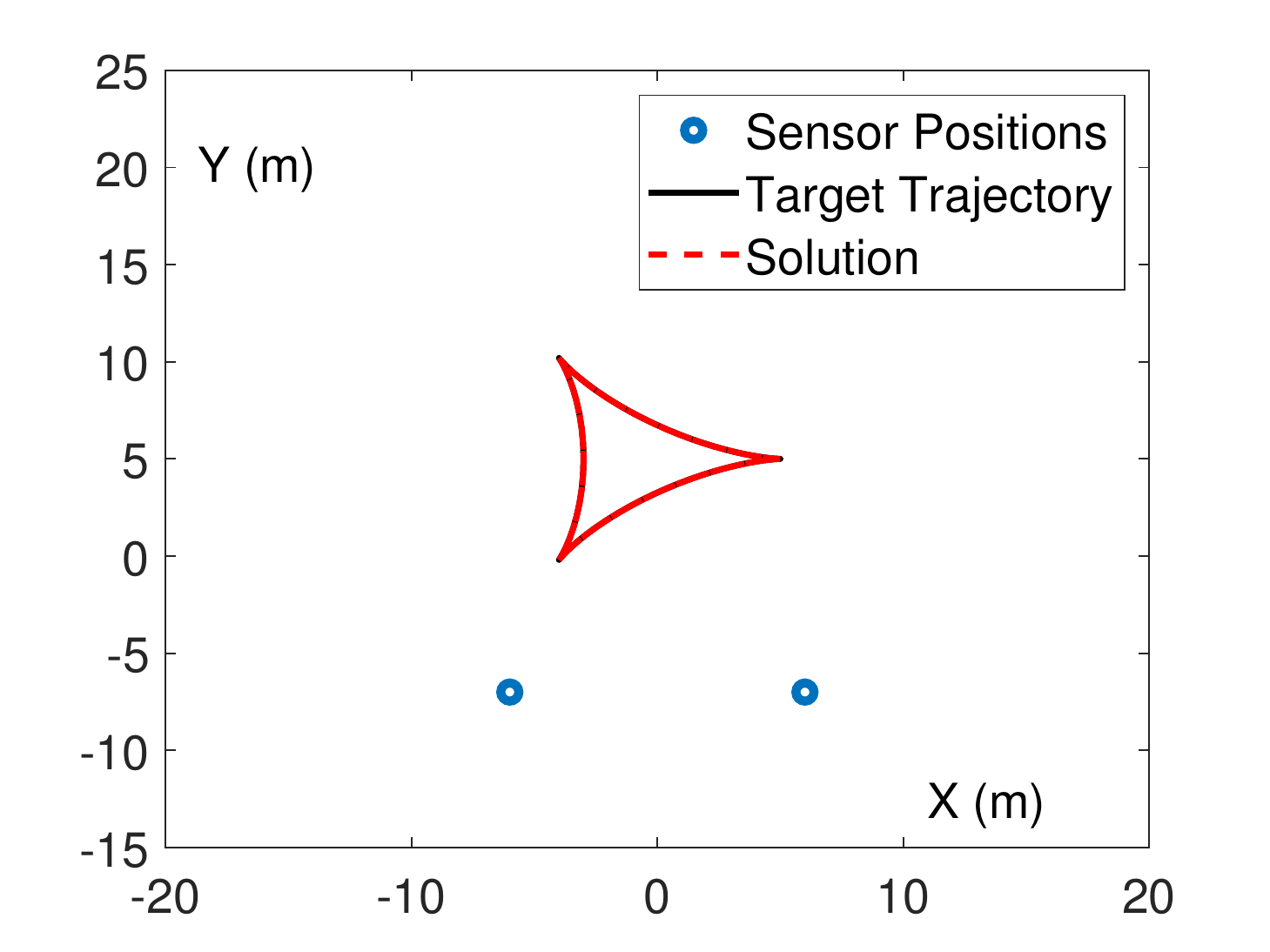}}
\subfigure[Position error of the AZTND model (\ref{RACZNNforAoA})]{\includegraphics[scale=0.33]{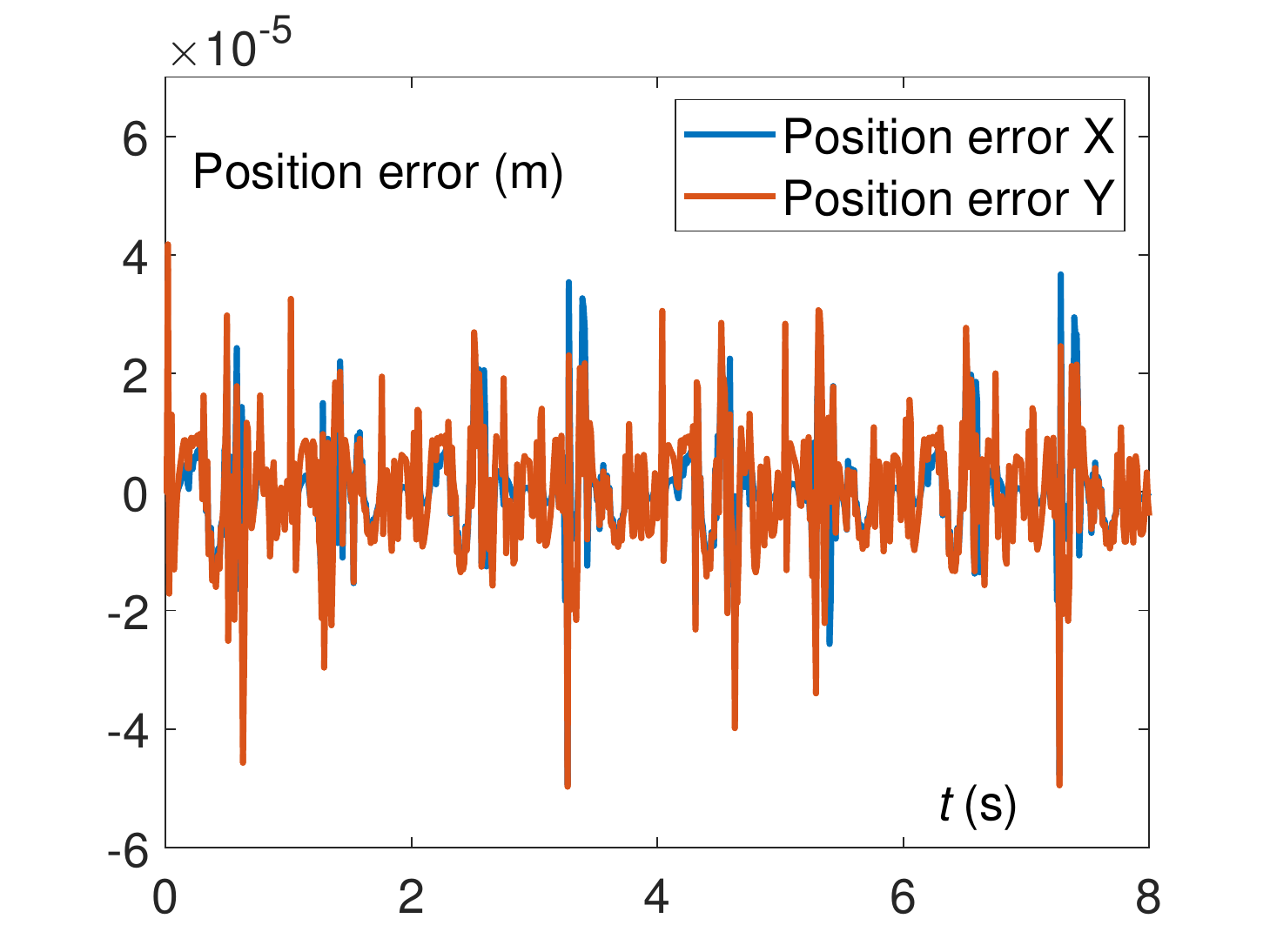}}
\subfigure[Position error of the OZNN model (\ref{OZNNforAoA})]{\includegraphics[scale=0.33]{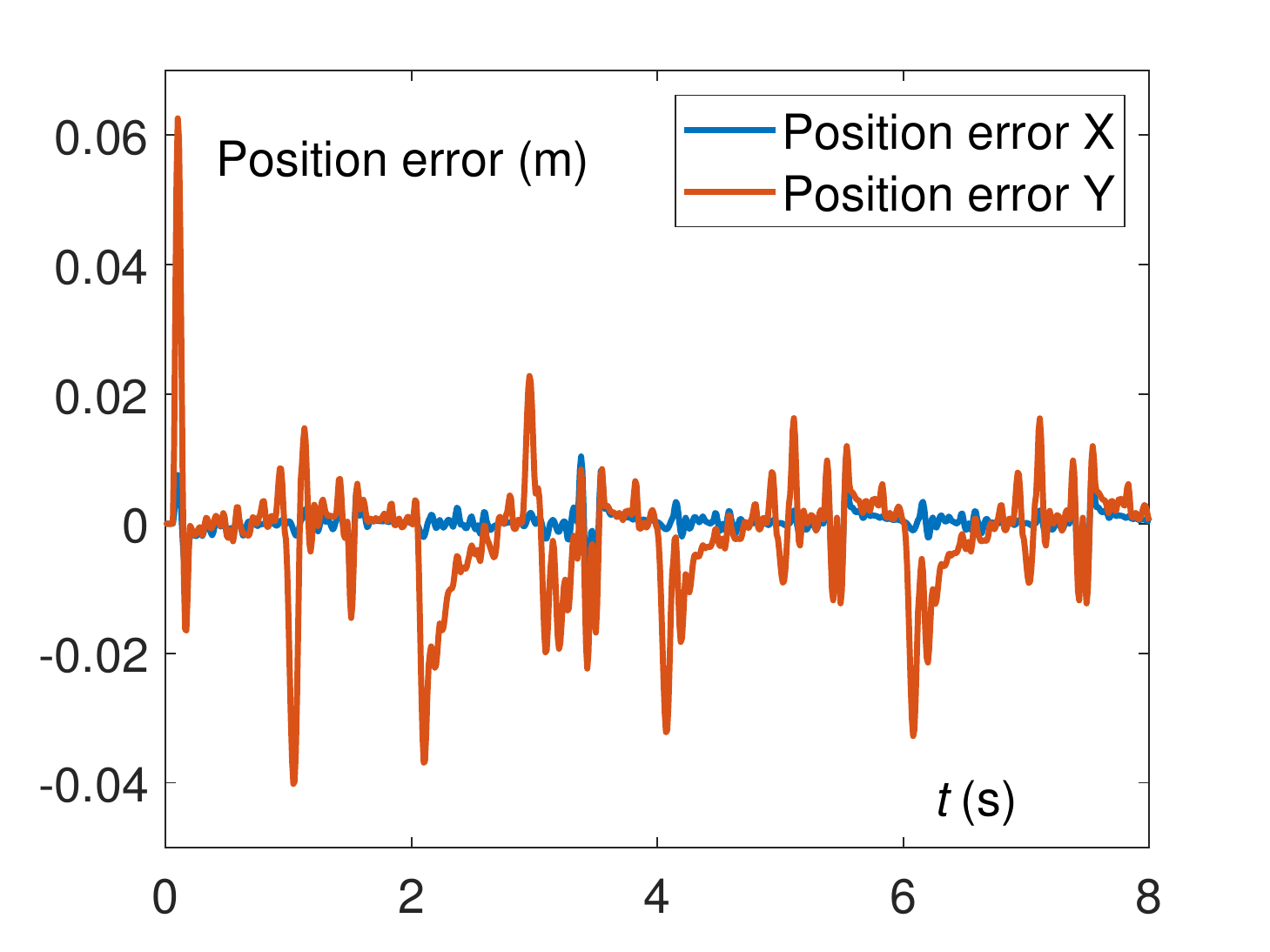}}
\subfigure[Theoretical and estimated trajectory]{\includegraphics[scale=0.33]{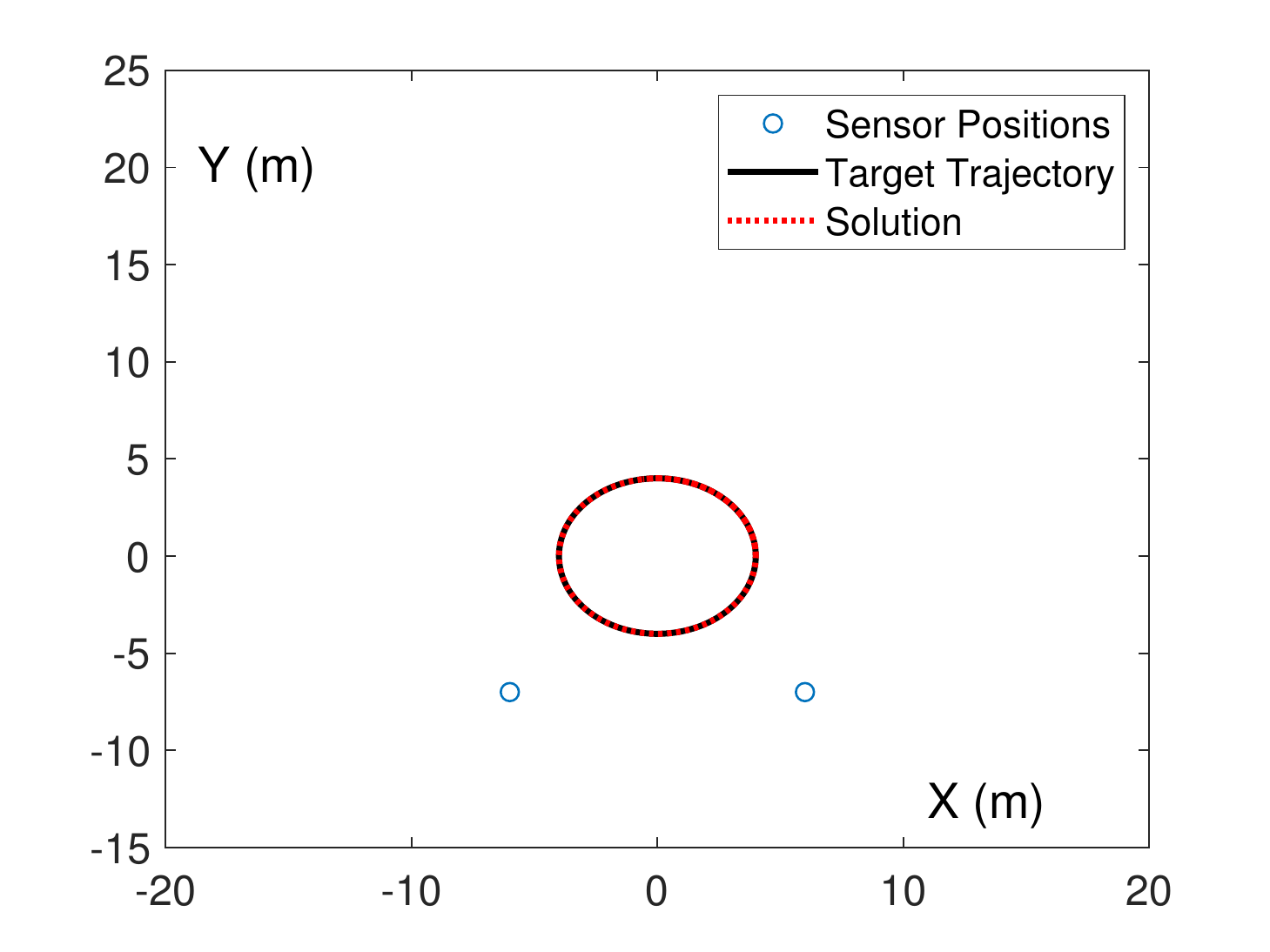}}
\subfigure[Position error of the AZTND model (\ref{RACZNNforAoA})]{\includegraphics[scale=0.33]{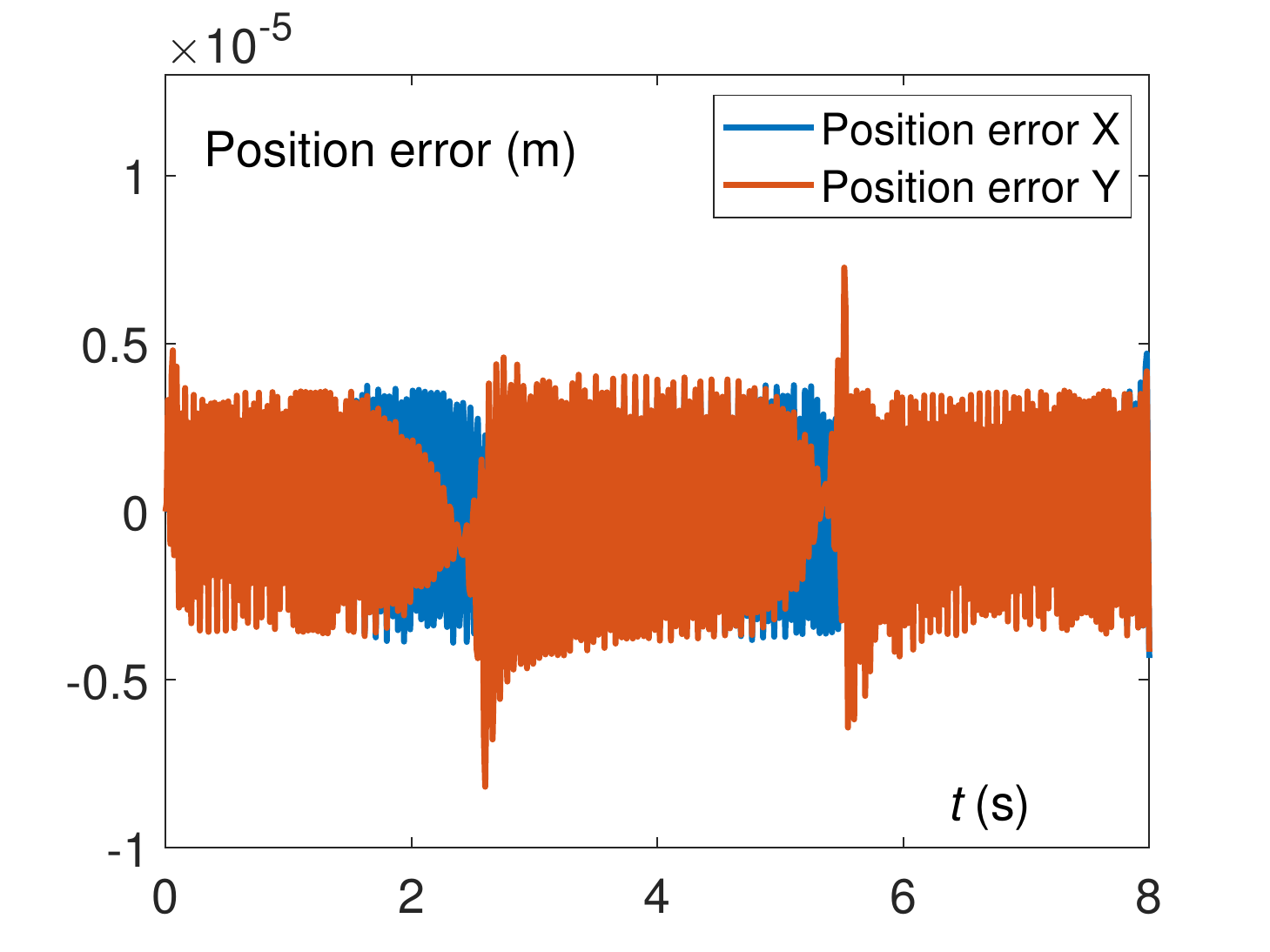}}
\subfigure[Position error of the OZNN model (\ref{OZNNforAoA})]{\includegraphics[scale=0.33]{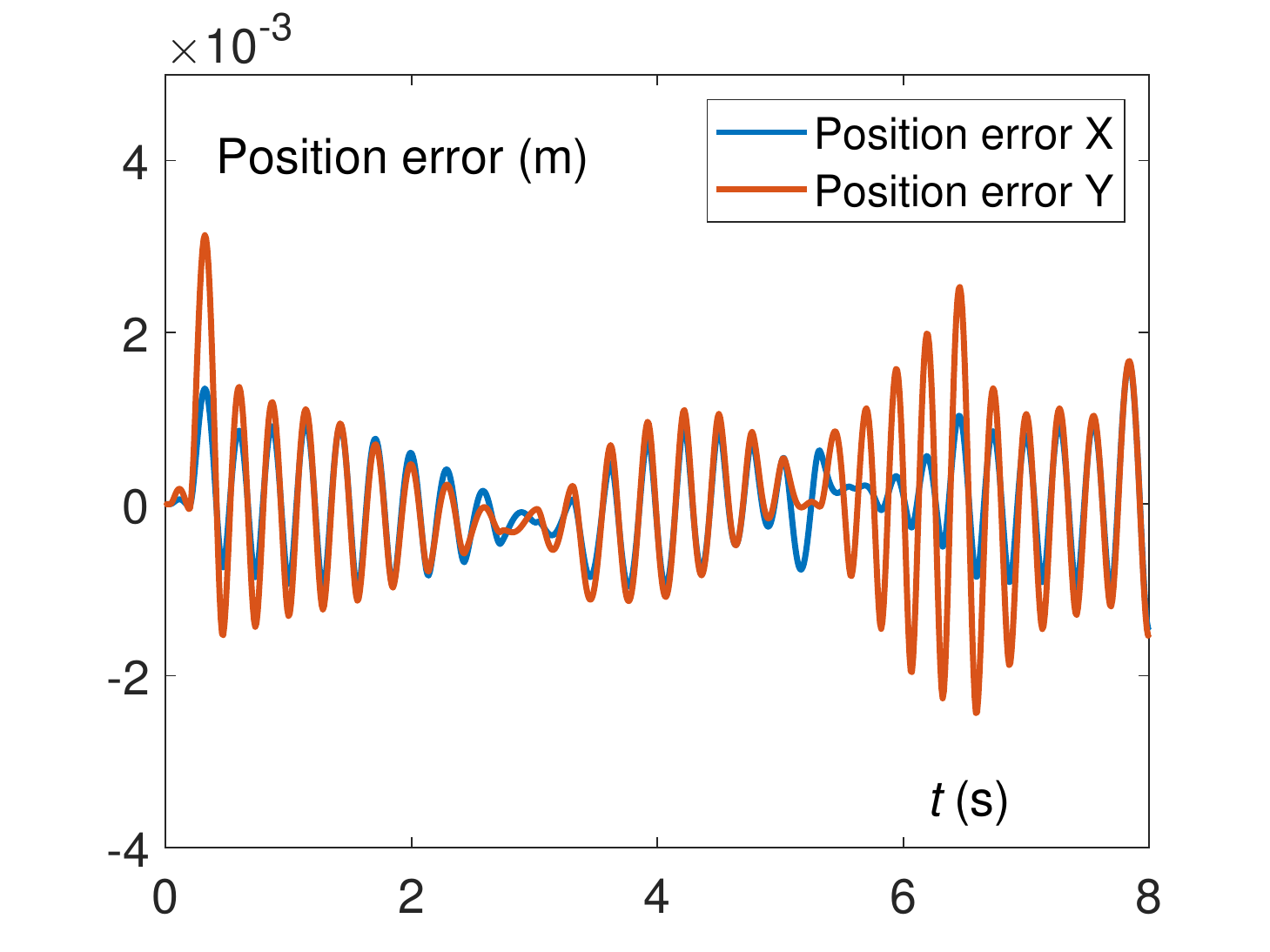}}
\caption{Visualization results of the AoA dynamic positioning scheme. (a) and (d) denoting solution result trajectories (red dashed line) generated by the AZTND model (\ref{RACZNNforAoA}) and the theoretical trajectory (solid black line) of the dynamic target. (b) and (e) representing the corresponding position error synthesized by the AZTND model (\ref{RACZNNforAoA}). (c) and (f) denoting the corresponding position error synthesized by the OZNN model (\ref{OZNNforAoA}).
}
\label{AoAResult}
\end{figure*}

\begin{figure}[htbp]\centering
\includegraphics[scale=0.62]{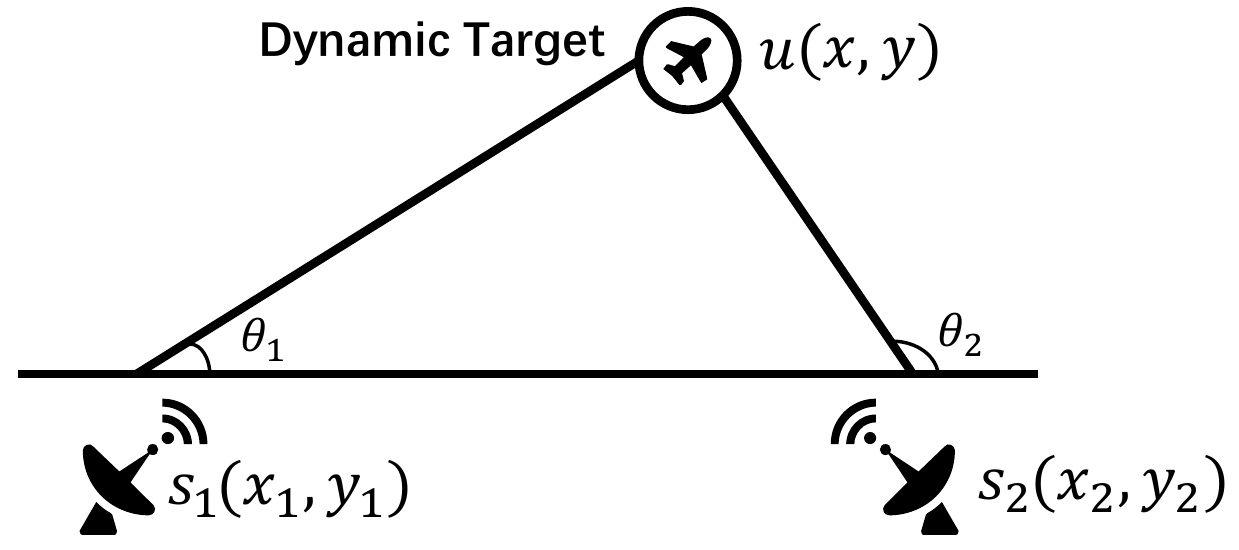}
\caption{
Schematic diagram of AoA dynamic positioning scheme.
}
\label{AoA}
\end{figure}

\subsubsection{AZTND Model with Noise}
The quantitative experiment visualization results synthesized by the AZTND model (\ref{RACZNN}) for solving the TVQM problem example (\ref{EA}) in the constant noise, linear noise, and bounded random noise cases are arranged in Figure \ref{ConNorm}, Figure \ref{LinearNorm}, and Figure \ref{RandNorm}, respectively. Besides, the adaptive scale coefficient and adaptive feedback coefficient of the AZTND model (\ref{RACZNN}) are $\xi(\vec{\epsilon}(t)) = ||\vec{\epsilon}(t)||_{\text{2}}^{3} + 5$ and $\kappa(\vec{\epsilon}(t)) = 5^{||\int_0^t\vec{\epsilon}(\delta)\text{d}\delta||_{\text{F}}}+5$, respectively. Based on three different situations, the following three situations will be analyzed and discussed in detail. Noticeably, the measurement noises perturbed by the AZTND model can be expressed as equation (\ref{RACZNNNoise}).
Firstly, the amplitude of the constant noise in example (\ref{EA}) is provided as $\vec{\vartheta}(t)=\vec{\vartheta}=[5,5]^{\text{T}}$. As shown in Figure \ref{ConNorm} (a), beginning with a random-generated initial value, even though the AZTND model (\ref{RACZNN}) is disturbed by the constant noise, its system residual error $||\vec{\epsilon}(t)||_{\text{2}}$ still accurately converge to the theoretical solution. Meanwhile, Figure (\ref{ConNorm}) (b) depicts that the solution accuracy of the GNN model (\ref{GNNCompare}), PTCZNN model (\ref{PTCZNNCompare}), and NCZNN model (\ref{NCZNNCompare}) remain at a relatively high level. Secondly, the quantitative experimental simulation results of the AZTND model (\ref{RACZNN}) solving the TVQM problem example (\ref{EA}) under linear noise $\vec{\vartheta}(t)=\vec{\vartheta}(t)\in \mathbb{R}^{n}$ interference are arranged in Figure \ref{LinearNorm}. Suffering from the interference of linear noise $\vec{\vartheta}(t)$, the proposed AZTND model (\ref{RACZNN}) not only has the highest solution accuracy compared with the other four comparative models, but also its convergent speed is the second only to the GNN model (\ref{GNNCompare}). Finally, the quantitative experiment simulation results of the AZTND model (\ref{RACZNN}) solving the TVQM problem example (\ref{EA}) under bounded random noise $\vec{\vartheta}(t)=\vec{\varrho}(t)\in [0.5, 3]^2$ interference are arranged in Figure \ref{LinearNorm}. As displayed in Figure \ref{LinearNorm} (a) and (b), although the proposed AZTND model (\ref{RACZNN}) is interrupted by bounded random noise $\vec{\varrho}(t)$, its average steady-state residual error still maintains a high accuracy, specifically, of order $10^{-2}$. By contrast, other commonly used neural network models, $i.e.$, the GNN model (\ref{GNNCompare}), PTCZNN model (\ref{PTCZNNCompare}), and NCZNN model (\ref{NCZNNCompare}) remain relatively large residual errors, specifically, of order $10^{-1}$. Furthermore, a similar conclusion is drawn from Figure \ref{RandNorm}. Therefore, the conclusion is drawn that the proposed AZTND model (\ref{RACZNN}) has higher robustness and stability when facing different noises than other state-of-the-art neural network models.

\begin{figure*}[htbp]\centering
\subfigure[]{\includegraphics[scale=0.4]{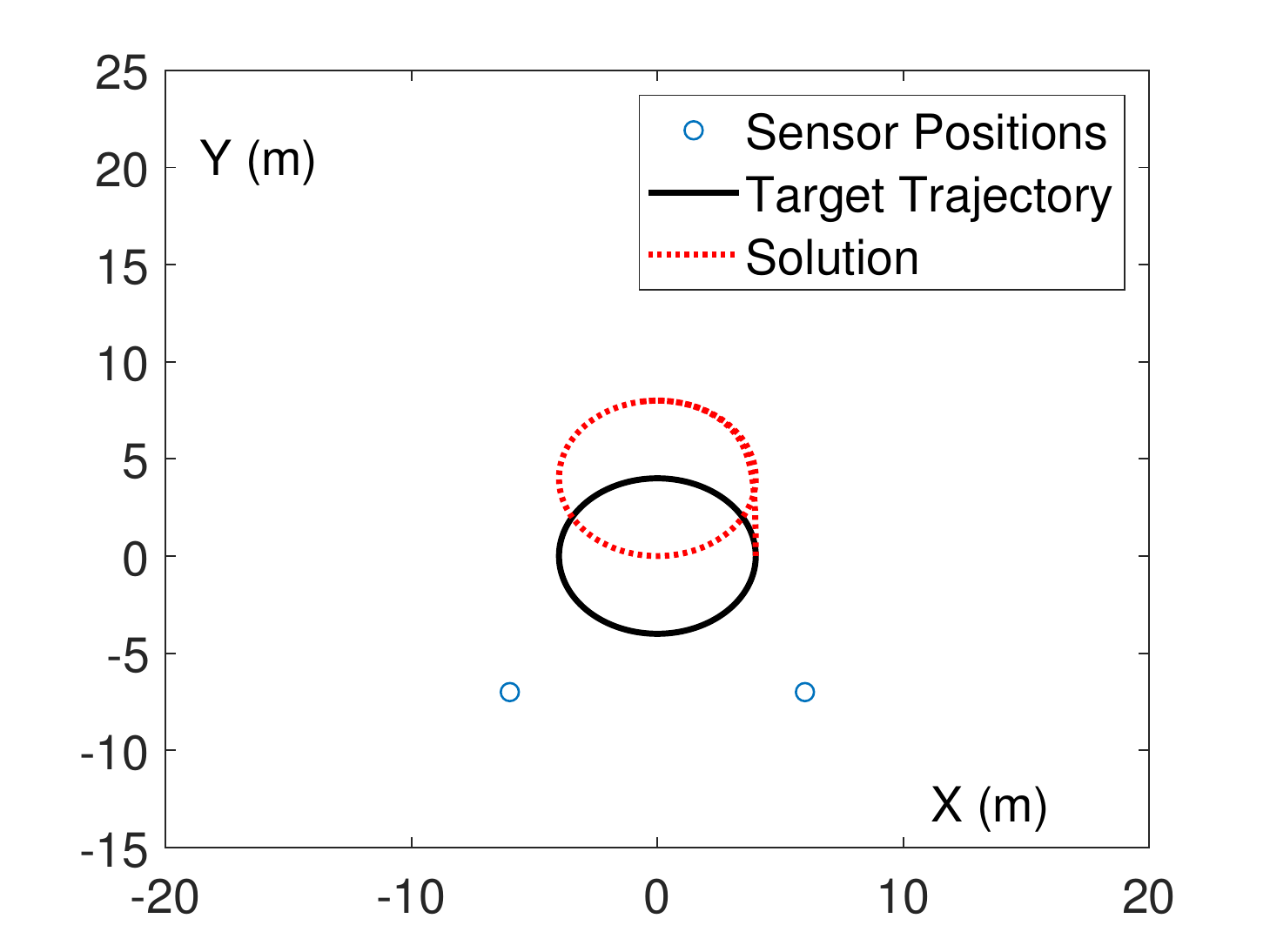}}
\subfigure[]{\includegraphics[scale=0.4]{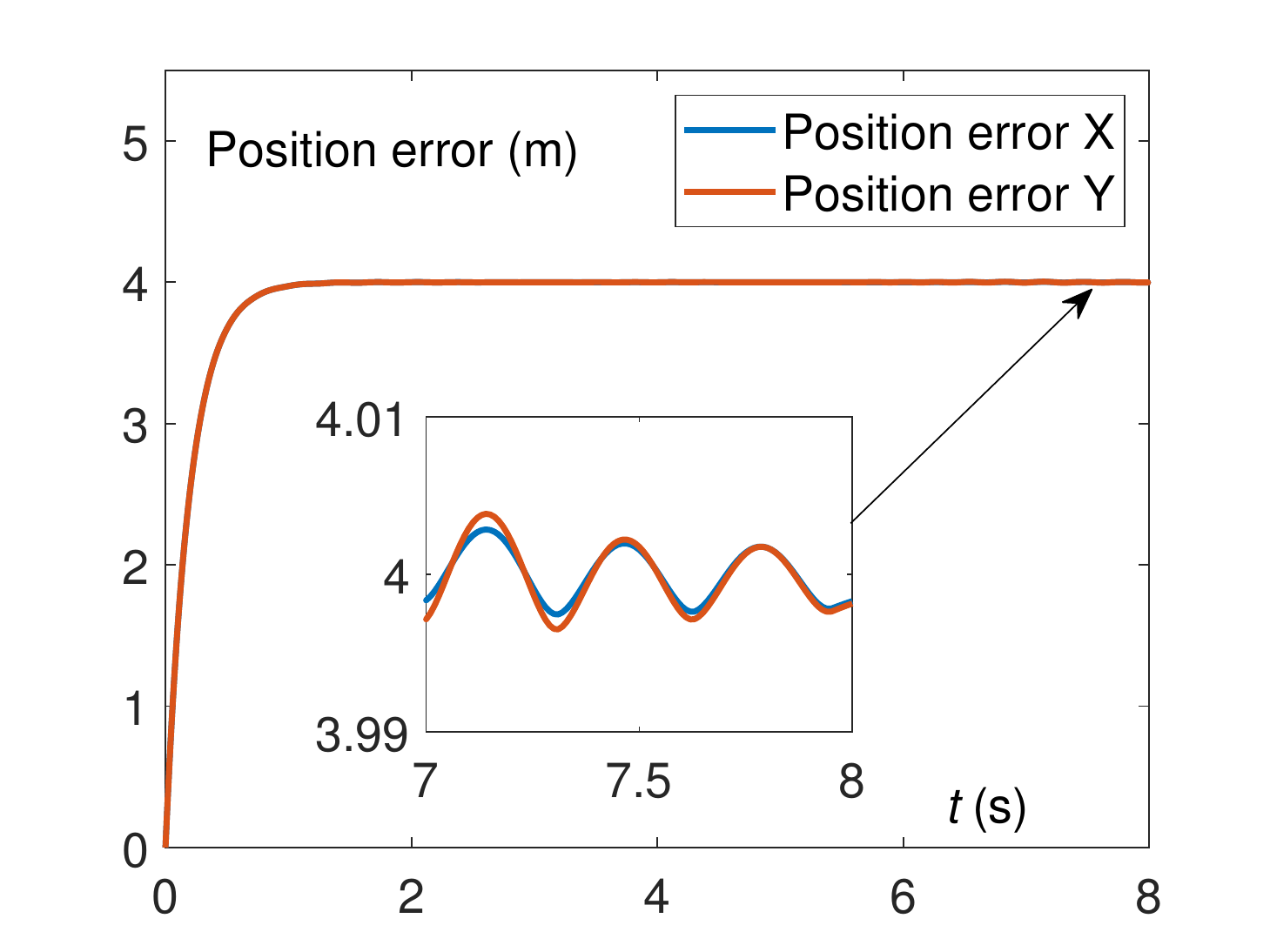}}
\subfigure[]{\includegraphics[scale=0.4]{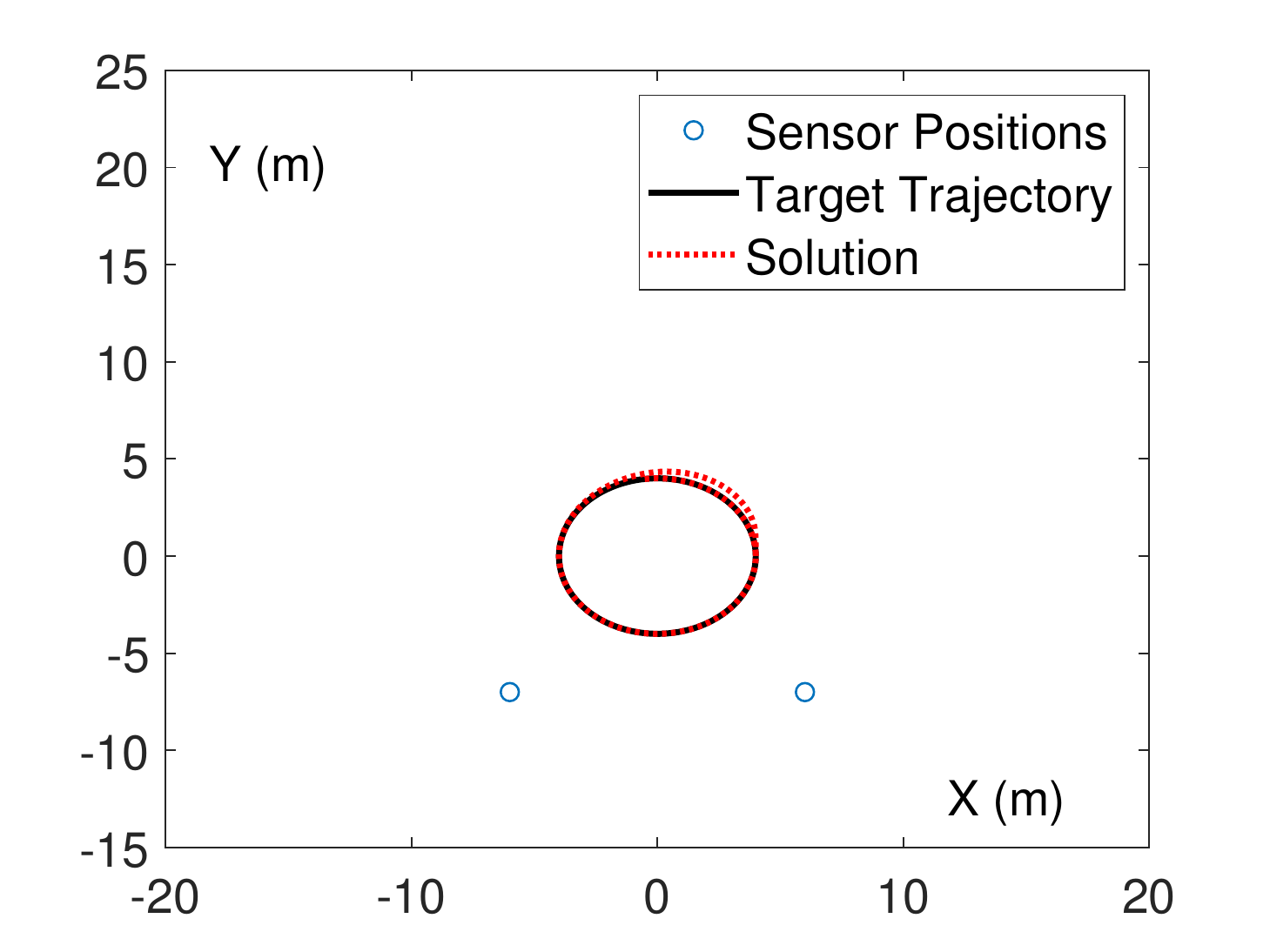}}
\subfigure[]{\includegraphics[scale=0.4]{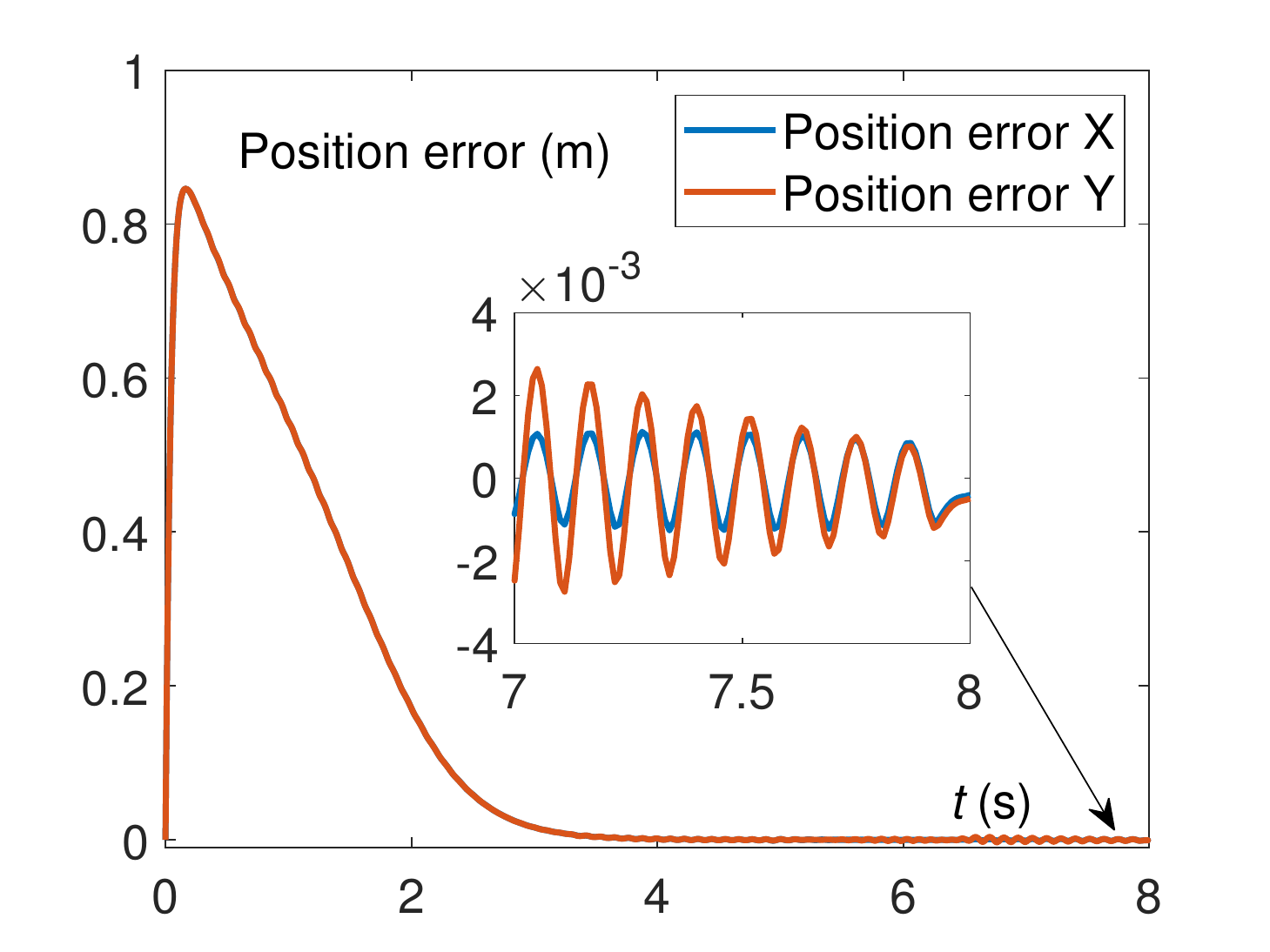}}
\caption{Visualization results of the AoA dynamic positioning scheme with constant noise $\vec{\vartheta}(t)=[20,20]^{\text{T}}$ disturbance. (a) and (c) showing the solution result trajectories generated by the AZTND model (\ref{RACZNNforAoA}) and the OZNN model (\ref{OZNNforAoA}). (b) and (d) denoting the corresponding position errors of the OZNN model (\ref{OZNNforAoA}) and the AZTND model (\ref{RACZNNforAoA}).
}
\label{AoANoiseResult}
\end{figure*}

\subsection{AZTND Applied to Target Tracking Scheme}
In this part, an angle-of-arrival (AoA) target tracking scheme based on the AZTND model (\ref{RACZNN}) is presented, which is widely used in navigation, guidance, and localization system. The following two-dimensional (2-D) AoA positioning principle is detailed to construct the target tracking scheme further. As shown in Figure \ref{AoA}, rays emitted by the observation base stations $s_1$ and $s_2$ will pass through the dynamic target $u$, and the intersection of two rays is the dynamic target's position. The position coordinate $u(x,y)$ of the dynamic target can be solved by calculating the angles of arrival $\theta_1$ and $\theta_2$  observation base stations to the dynamic target.

Consequently, the geometric relationship between the angle of arrival between the observation base station and the target is expressed as
\begin{equation}\label{AoATan}
\text{tan}\theta_i(t) = \frac{y(t)-y_i}{x(t)-x_i}.
\end{equation}
Noting that parameters $\theta_i(t)$ and $(x(t),y(t))$ represent real-time changing angles and coordinates. Expand and rearrange the equation (\ref{AoATan}) to get $y_i-x_i\text{tan}\theta_i(t)=-x(t)\text{tan}\theta_i(t)+y(t)$. Therefore, the following linear equation is obtained:
\begin{eqnarray}\label{AOAOriginal}
    \begin{bmatrix}
        -\text{tan}(\theta_1(t))& 1\\
        -\text{tan}(\theta_2(t))& 1\\
        \vdots &\vdots\\
        -\text{tan}(\theta_n(t))& 1
    \end{bmatrix}
    \begin{bmatrix}
        x(t)\\
        y(t)
    \end{bmatrix}=
    \begin{bmatrix}
        y_1-x_1\text{tan}(\theta_1(t))\\
        y_2-x_2\text{tan}(\theta_2(t))\\
        \vdots\\
        y_n-x_n\text{tan}(\theta_2(t))
    \end{bmatrix},
\end{eqnarray}
where the parameter $n$ denotes the number of observation base stations, and we further express equation (\ref{AOAOriginal}) as the following equation:
\begin{equation}\label{AoALinear}
F(t)\vec{g}(t)=\vec{h}(t),
\end{equation}
where the time-varying matrix $F(t)\in\mathbb{R}^{n\times 2}$, time-varying vector $\vec{g}(t)=[x(t),y(t)]^{\text{T}}$ and $\vec{h}(t)\in\mathbb{R}^{n}$. Besides, the error function for equation (\ref{AoALinear}) is written as $\vec{\epsilon}(t)=F(t)\vec{g}(t)-\vec{h}(t)$. Subsequently, the AZTND model (\ref{RACZNN}) for the target tracking scheme is presented as
\begin{eqnarray}\label{RACZNNforAoA}
	\begin{split}
		 F(t)\vec{\dot g}(t)=&-\dot F(t)\vec{g}(t)-\vec{\dot h}(t)-\xi(\vec{\epsilon}(t))\big{(}F(t)\vec{g}(t)\\
		&-\vec{h}(t)\big{)} -\kappa(\vec{\epsilon}(t))\int_{0}^{t}(F(\delta)\vec{g}(\delta)+\vec{h}(\delta))\text{d}\delta.
	\end{split}
\end{eqnarray}
For comparison, the OZNN for the AoA target tracking scheme is introduced as follows:
\begin{equation}\label{OZNNforAoA}
	F(t)\vec{\dot g}(t)=-\dot F(t)\vec{g}(t)-\vec{\dot h}(t)-\gamma\big{(}F(t)\vec{g}(t)-\vec{h}(t)\big{)},
\end{equation}
where the parameter $\gamma$ represents the scale parameter. The corresponding simulation results are provided in Figure \ref{AoAResult} and Figure \ref{AoANoiseResult}. Figure \ref{AoAResult} shows the target trajectory and system position error obtained by the target tracking scheme based on the OZNN model and the proposed AZTND model. Starting from the initial coordinate point $u_0$, the position error generated by the AZTND model (\ref{RACZNNforAoA}) converges to the $10^{-5}$ order, which is more accurate than the $10^{-3}$ order obtained by the OZNN model (\ref{OZNNforAoA}). Besides, as described in Figure \ref{AoANoiseResult} (a) and (c), the trajectory obtained by the OZNN model (\ref{OZNNforAoA}) does not converge to the true dynamic target trajectory when disturbed by constant noise $\vec{\vartheta}(t)=[20,20]^{\text{T}}$, while the trajectory generated by the AZTND model (\ref{RACZNNforAoA}) is well consistent with the dynamic target trajectory. Figure \ref{AoANoiseResult} (b) and (d) show that the AZTND model (\ref{RACZNNforAoA}) with constant noise $\vec{\vartheta}(t)$ interference still converges to order $10^{-3}$, while the OZNN model (\ref{OZNNforAoA}) diverges. In general, this part fully demonstrates that the AZTND model (\ref{RACZNNforAoA}) is effectively applied to the AoA target tracking scheme regardless of whether there is noise or not.

\section{Conclusions}\label{Conclusion}
The adaptive zeroing-type neural dynamics (AZTND) model has been proposed in this paper to solve the time-varying quadratic minimization (TVQM) problem in a perturbed environment. Unlike the original zeroing neural network models, the scale coefficient and feedback coefficient of the AZTND model has been presented from the perspective of adaptive optimization to expedite the global convergence and enhance the robustness of the model. Furthermore, this paper has presented the corresponding theorem and proof procedures from the stability perspective to investigate the global convergence of the AZTND model. Then, the corresponding numerical experiment is designed and executed, and the numerical results and visualization results of the investigation are given in tables and figures, respectively. Finally, the potential of the AZTND model in practical applications has been shown, and the simulative experiment demonstrates the effectiveness and superiority of the target tracking scheme based on the AZTND model.



\begin{thebibliography}{99}
\bibitem{AppEnergy}
Killian M., Zauner M., Kozek M., 2018. Comprehensive smart home energy management system using mixed-integer quadratic programming. Applied Energy 222, 662--672.


\bibitem{NCZNNCompare}
Jiang C., Xiao X., Liu D., Huang H., Xiao H., Lu H., 2021. Nonconvex and bound constraint zeroing neural network for solving time-varying complex-valued quadratic programming problem. IEEE Trans. Ind. Informat. 17, 6864--6874.

\bibitem{AppRobort}
Zhang Z., Li Z., Zhang Y., Luo Y., Li Y., 2015. Neural-dynamic-method-based dual-arm CMG scheme with time-varying constraints applied to humanoid robots. IEEE Trans. Neural. Netw. Learn. Syst. 26, 3251--3262.

\bibitem{VPZNNCompare}
Zhang Z., Lu Y., Zheng L., Li S., Yu Z., Li Y., 2018. A new varying-parameter convergent-differential neural-network for solving time-varying convex QP problem constrained by linear-equality. IEEE Trans. Autom. Control. 63, 4110--4125.


\bibitem{AppMessagePass}
Ruozzi N., Tatikonda S., 2013. Message-passing algorithms for quadratic minimization. J. Mach. Learn. Res. 14, 2287--2314.

\bibitem{XXCOne}
Xiao X., Xiong N. N., Lai J., Wang C. D., Sun Z., Yan J., 2021. A local consensus index scheme for random-valued impulse noise detection systems. IEEE Trans. Syst. Man. Cybern. Syst. 51, 3412--3428.

\bibitem{ZNNProposed}
Zhang Y., Mu B., Zheng H., 2013. Link between and comparison and combination of Zhang neural network and quasi-Newton BFGS method for time-varying quadratic minimization. IEEE Trans. Cybern. 43, 490--503.

\bibitem{XJ}
Xiao X., Jiang C., Lu H., Jin L., Liu D., Huang H., Pan Y., 2020. A parallel computing method based on zeroing neural networks for timevarying complex-valued matrix Moore-Penrose inversion. Inf. Sci. 524, 216--228.

\bibitem{LiuMei}
Lui M., Chen L., Du X., Jin L., Shang M., 2021. Activated gradients for deep neural networks. IEEE Trans. Neural Netw. Learn. Syst. http//doi.org/10.1109/TNNLS.2021.3106044.

\bibitem{Qi}
Qi Y., Jin L., Wang Y., Xiao L., Zhang J., 2020. Complex-valued discrete-time neural dynamics for perturbed time-dependent complex quadratic programming with applications. IEEE Trans. Neural Netw. Learn. Syst. 31, 3555--3569.

\bibitem{Wei}
Wei L., Jin L., Yang C., Chen K., Li W., 2021. New noise-tolerant neural algorithms for future dynamic nonlinear optimization with estimation on hessian matrix inversion. IEEE Trans. Syst. Man Cybern. Syst. 51, 2611--2623.


\bibitem{AppWideUse}
Xiao L., Li S., Yang J., Zhang Z., 2018. A new recurrent neural network with noise-tolerance and finite-time convergence for dynamic quadratic minimization. Neurocomputing 285, 125--132.

\bibitem{FTZNNXiao}
Xiao L., Tan H., Jia L., Dai J., Zhang Y., 2020. New error function designs for finite-time ZNN models with application to dynamic matrix inversion. Neurocomputing 402, 395--408.

\bibitem{PTCZNNCompare}
Li W., Ma X., Luo J., Jin L., 2021. A strictly predefined-time convergent neural solution to equality- and inequality-constrained time-variant quadratic programming. IEEE Trans. Syst. Man Cybern. Syst. 51, 4028--4039.

\bibitem{Huang}
Miao P., Wu D., Shen Y., Zhang Z., 2019. Discrete-time neural network with two classes of bias noises for solving time-variant matrix inversion and application to robot tracking. Neural Comput. Applic. 31, 4879--4890.

\bibitem{MZNNCompare}
Jin L., Zhang Y., Li S., Zhang Y., 2016. Modified ZNN for time-varying quadratic programming with inherent tolerance to noises and its application to kinematic redundancy resolution of robot manipulators. IEEE Trans. Ind. Electron. 63, 6978--6988.

\bibitem{Alawad}
Alawad M., Lin M., 2019. Survey of stochastic-based computation paradigms. IEEE Trans. Emerg. Topics. Comput. 7, 98--114.

\bibitem{DeepRL}
He K., Zhang X., Ren S., Sun J., 2016. Deep residual learning for image recognition. IEEE
Conf. Comp. Vis. Patt. Recogn. 770--778.

\bibitem{IterCompareOne}
Liu Q., Wang J., 2015. A projection neural network for constrained quadratic minimax optimization.IEEE Trans. Neural Netw. Learn. Syst. 26, 2891--2900.

\bibitem{ACGNNCompare}
Liao S., Liu J., Xiao X., Fu D., Wang G., Jin L., 2020. Modified gradient neural networks for solving the time-varying Sylvester equation with adaptive coefficients and elimination of matrix inversion. Neurocomputing 379, 1--11.

\bibitem{OZNNCompare}
Xiao X., Fu D., Wang G., Liao S., Qi Y., Huang H., Jin L., 2020. Two neural dynamics approaches for computing system of time-varying nonlinear equations. Neurocomputing 394, 84--94.

\bibitem{GNNCompare}
Chen Y., Yi C., Qiao D., 2013. Improved neural solution for the Lyapunov matrix equation based on gradient search. Inf. Process Lett. 13, 876--881.
\end{thebibliography}
\end{document}